\documentclass[amsmath,amssymb]{article}
\usepackage[numbers,sort&compress]{natbib}
\usepackage{amsfonts,latexsym,graphicx}

\addtocounter{secnumdepth}{1}

\usepackage{epstopdf}
\usepackage[bookmarksnumbered,colorlinks,bookmarks,citecolor=blue,linkcolor=blue,breaklinks,linktocpage]{hyperref}
\usepackage[all]{hypcap}
\usepackage{subfigure}

\usepackage{dcolumn}
\usepackage{bm}
\usepackage{amsfonts,latexsym,graphicx}
\usepackage[applemac]{inputenc}
 \usepackage[polutonikogreek,english]{babel}
\usepackage{color}
\usepackage{multirow}
\usepackage{subfigure}

\usepackage{amsthm}

\usepackage{epigraph}

\usepackage{diagrams}
\usepackage[all]{xy}

\usepackage{tikz-cd}
\usetikzlibrary{arrows}
\tikzset{commutative diagrams/.cd,arrow style=tikz,diagrams={>=latex'}}
\usepackage{graphicx}

\theoremstyle{plain}
\newtheorem{theorem}{Theorem}[section]

\newtheorem{conjecture}[theorem]{Conjecture} 

\theoremstyle{definition}
\newtheorem{definition}[theorem]{Definition}

\theoremstyle{remark}

\newtheorem{example}[theorem]{Example}

\begin{document}
\title{Introduction to Gestural Similarity in Music. \\ An Application of Category Theory to the Orchestra}

\author{
Maria Mannone\footnote{manno012@umn.edu}
\\School of Music,
\\University of Minnesota,
Minneapolis, USA
\\ THIS IS NOT the latest version.
\\The latest version has been published
\\ by the Journal of Mathematics and Music.
}

\maketitle

\begin{abstract}
Mathematics, and more generally computational sciences, intervene in several aspects of music. Mathematics describes the acoustics of the sounds giving formal tools to physics, and the matter of music itself in terms of compositional structures and strategies. Mathematics can also be applied to the entire {\em making} of music, from the score to the performance, connecting compositional structures to acoustical reality of sounds. Moreover, the precise concept of {\em gesture} has a decisive role in understanding musical performance. In this paper, we apply some concepts of category theory to compare gestures of orchestral musicians, and to investigate the relationship between orchestra and conductor, as well as between listeners and conductor/orchestra. To this aim, we will introduce the concept of {\em gestural similarity}. The mathematical tools used can be applied to gesture classification, and to interdisciplinary comparisons between music and visual arts.
\end{abstract}


\section*{Introduction}

The topic of {\em musical gestures} in performance and composition is the object of an increasing interest from scholars \citep{godoy_leman, visi_miranda_2, greek_2, ToM_III}. Interest in gestural comparison is a natural consequence. New research about the mathematical description of gestures starts from piano, moves to orchestral conducting, and reaches singing as `inner movement' \citep[Chapter 37]{ToM_III}. The topic of {\em musical similarity} is another important object of research \citep{cont, mannone_mus_nonMarkov, new_book_sim}. 

In music, gestures are the intermediary between musical thinking and acoustics.
Musical scores contain playing instructions and performers interact with their instruments via gestures to produce a specific sound result.
Gestures are also involved in musical perception \citep{lipscomb}.

In this article, we focus on the comparison, in terms of {\em gestural similarity}, among gestures of the orchestral musicians, conductor's gestures, and specialization of the conductor's instructions into the movements of each musician, using the framework of category theory \citep{macLane}. We will also briefly refer to visual arts, using the concept of gesture to connect them with music.


The interest of scholars in categories is motivated by their recent applications in physics and other sciences \citep{spivak}, and, more generally, in philosophy and diagrammatic thinking \citep{alunni}. Mathematical music theory that uses category theory has already been developed \citep{gm:TOM}, as well as differential calculus for musical analysis and musical gestures \citep{gm:TOM, mcm15, global}. In conclusion, we refer to fuzzy logic for the definition of the degree of gestural similarity.  The fuzzy logic had already been applied to music while investigating the relations between emotions and music \citep{friberg}.

\subsection*{The state of the art of the mathematical theory of musical gestures}



\begin{figure}
\centerline{
\includegraphics[width=9cm]{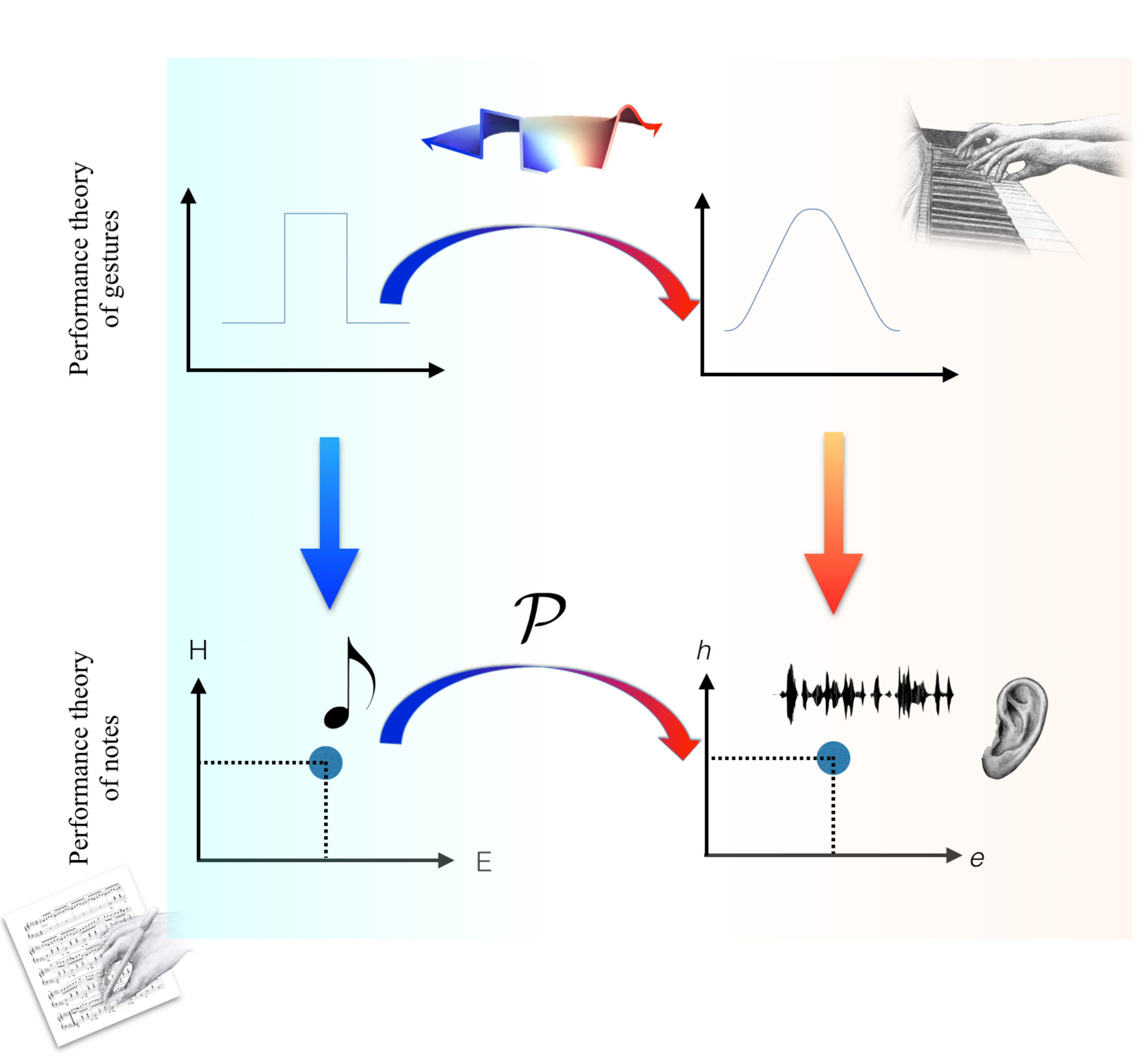}}
\caption{\label{scheme} \footnotesize The general scheme of mathematical performance theory, from gestures (top) to notes (bottom). From the original German notation adopted in \citep{gm:TOM}, $H$ indicates symbolic pitch, $h$ physical pitch, $E$ symbolic onset, $e$ physical onset. The other parameters are here omitted for graphical reasons. The graphs in the upper part of the image represent the gesture of the fingertip on piano keyboard, in the simplified case onset (in abscissa) and vertical position (in ordinate).
}
\end{figure}

What is a musical {\em gesture}? We can intuitively think of the movements of a dancer (gestures) while touching the floor at discrete points (notes), see Figure \ref{scheme}.

We can distinguish two different types of gestures: {\em symbolic} and {\em physical} gestures. The first are derived from the information contained in the score, and the second are the performer's real movements, see the top of Figure \ref{scheme}. Let us refer to the case of the pianist as a starting point. The symbolic gesture is derived from information such as the MIDI-like command ``play this key, at this time, and with this loudness.'' The symbolic gesture contains straight lines and, in principle, it can require an infinite speed for the transition from pressed to released key, and from one note to another. The physical gesture is given by the real, physically-possible movements made by performers in non-zero time and finite speed, represented by the hands' smooth paths (their smooth curves) in space and time.

Symbolic gestures can be transformed into physical ones via an ideal connecting surface, the {\em world-sheet}. This formal tool, coming from string theory \citep{zwiebach} in theoretical physics---in fact, we are comparing the continuous curves of gestures to strings---is interpreted here as a {\em hypergesture}.\footnote{In string theory, point-like particles are substituted by vibrating strings. In \citep{mcm15} the paradigm of string theory is used to describe musical performance in terms of gestures, rather than notes and sounds as isolated events. The `particles' correspond to the notes. To understand music as unfolding in time, rather than as isolated ``points in time,'' we can study gestures to understand music. However, the ``strings'' here are not vibrating: thus, such a reference to string theory is just a general metaphor.}
See Section \ref{hypergesture} for more details. However, because a gesture is a trajectory in a configuration space, we can formally study symbolic and physical gesture curves in the framework of Banach spaces and functional analysis. This would constitute a valid, and more precise, alternative approach to the string-theoretical one.

In this article, we refer more specifically to similarity between physical gestures. The concept of {\em similarity}, and even the term itself, had already been used in several musical contexts \citep{new_book_sim}, for example while referring to notes and motives: we can think of ``inner'' comparisons inside the same musical composition \citep{cont}, of motives' classification \citep{mazzola_motif} and compositions' comparisons, including an adaptation of the formalism of physics to study memory \citep{mannone_mus_nonMarkov, mannone_physics}.

Similarity can also be introduced for gestures, while dealing with {\em classification}.\footnote{For a precise but non-mathematical approach to gesture classification in the framework of mixed music (electroacoustic and with traditional instruments), see \citep{bachrata}; the term {\em similarity} is used to qualitatively describing imitational gestures, see \citep{semiotics2}.} How can we {\em compare} a gesture with other gestures of the same kind (piano with piano), or other kinds (piano with violin)? We may also wonder how to connect musical gestures with other forms of artistic expression (for example, drawing), or with external references, as in expressive gesture performance studies \citep{friberg}.

In this article we compare gestures having the same (or a similar) generator. We will make implicit use of homotopic transformation (see Paragraph \ref{hom} in the Appendix), as also proposed in contexts of gesture following \citep{greek_2, bevilacqua}. Let us consider the following example. When the orchestral conductor signals a {\em forte} for {\em tutti}, every musician makes a gesture or a combination of the more appropriate gestures depending on the technique of his or her musical instrument, to obtain acoustical spectra inside the desired range of loudness and timbre. Their gestures are thus {\em similar} (in loudness). We can intuitively find examples of such a concept also from composition. For example, while thinking of a delicate gesture, even outside music, the composer can write down a combination of notes, dynamics, tempi, to suggest to the performer a {\em delicate} final result. The performer then recreates, by using his or her knowledge of gestural technique, a sound as close as possible to the one {\em thought} by the composer\footnote{Of course, in the limit that such a thinking can be inferred from the score.} and filtered by the sensibility of the performer.

An inverse mechanism affects listening and music perception. A listener judges a musical performance as `expressive' if he or she {\em feels}, also beyond music, a gestural reference to be translated into mental images, a mental feeling. We can here refer to the concept of emotions and words viewed as hidden gestures \citep{wieland}, and refer to semiotics.\footnote{We will not go deeper into detail with semiotics. The passage from a simple, instrumental movement to an {\em expressive} movement may raise issues about semiotics. For example, a ``caressing'' piano touch is not only finalized to get a sound, but a sound with a specific, soft timbre, that also carries {\em a meaning}, the meaning of a caressing gesture. Visualization of such a gesture is also relevant for the perception of the intended meaning. An accurate study of this field should require a separate and detailed description, as well as some perception experiments, to substantiate the connection between our mathematical approach with the musicological studies in the field \citep{semiotics, semiotics2}.}
There can be cases where there is not any `human' gestural generator, as for electronic music, but the feeling of a clear gesture is recreated in the mind of the listener, and other cases when this process is not possible. Some studies compare `human' and electronic gestures, and describe analysis and creative developments of electronic music with respect to traditional musical instruments \citep{bachrata}. Even the {\em artistic inspiration} may be seen in light of gestural similarity, extending the language of artistic aesthetics to include the `ineffable' inside the `calculable' and understandable.

This article is structured as follows. After a list of preliminary mathematical concepts (Section \ref{list}), a first approach to gestural similarity is presented (Section \ref{orchestra gest sim}), with a discussion of the role of conductor (Section \ref{role conductor}), composer to conductor (Section \ref{adjoint}), and some short references to music and visual arts (Section \ref{visual}).
Some final remarks about the possible role of fuzzy logic in classification problems, as elements of future research, and an appendix with some comments on homotopic transformation and more mathematical details on musical gesture theory, conclude the article.

\section{Gestures and hypergestures}\label{list}

We need some preliminary mathematical tools before starting our analysis and investigation of gestural similarity.

\subsection{Mathematical definition of gesture}\label{mathematical definition of gesture}

\begin{figure}
 \centerline{
\includegraphics[width=10cm]{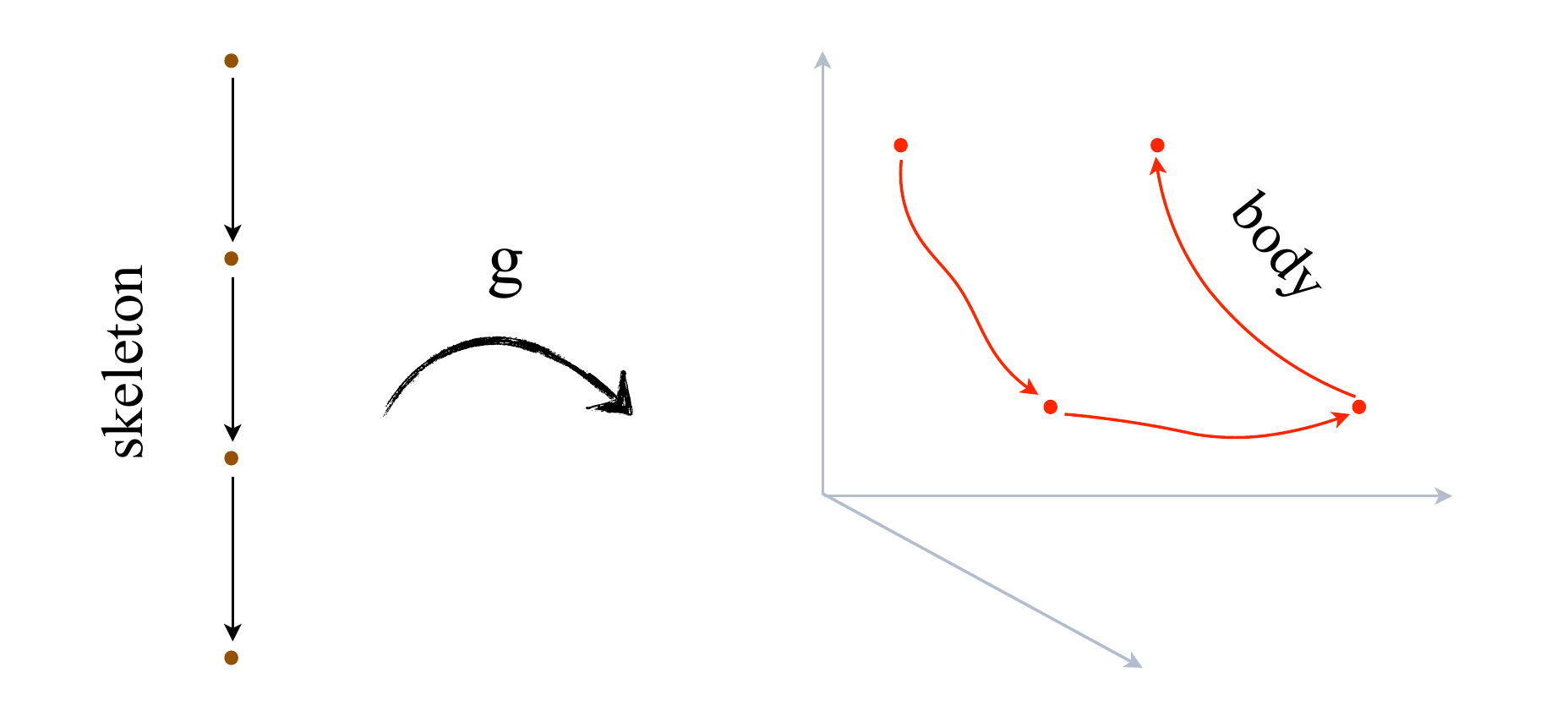}}
\caption{\label{gesture} \footnotesize An example of a gesture, a mapping from an oriented graph to a system of continuous curves in a topological space, from \citep{mazzola_andreatta}. We can compare this approach to a topic of neuroscience, about how to embed time-discrete symbolrefic processes into continuous time of neural processes \citep{graben2}.}
\end{figure}

Gestures have been mathematically defined as mappings from directed graphs to systems of continuous curves in topological spaces \citep{mazzola_andreatta}. More precisely, we start from an abstract system of points and connecting arrows which is a directed graph $\Delta$, as shown on the left of Figure \ref{gesture}. We  map this directed graph to a system of continuous curves\footnote{A curve $c$ in $X$ is a continuous function $c : I \rightarrow X$, where $I = [0, 1]$ is the real unit interval.} in a space\footnote{The space $\vec X$ is a topological space, e.g. spacetime.} $X$, with the same configuration as the directed graph $\Delta$, as shown on the right of Figure \ref{gesture}. The directed graph $\Delta$ is called the gesture’s {\em skeleton}, whereas the system of continuous curves constitutes the {\em body} of the gesture in $\vec X$. With $\vec X$ we indicate a set of mappings in $X$ (when we add a topology, $\vec X$ becomes the space of curves in $X$). A gesture is then a mapping $g:\Delta\rightarrow\vec X$.
We denote the space of all gestures from $\Delta$ to $\vec X$ as $\Delta@\vec X$, notation used in \citep{mazzola_andreatta} for $Hom(\Delta,\vec X)$. When we equip $\Delta@\vec X$ with a topology, we denote it as $\Delta\vec@ X$. This is the space of gestures.
A curve in the space of gestures (in $\Delta\vec@X$) is a gesture of gestures, and is called {\em hypergesture}.\footnote{This is the name introduced in \citep{mazzola_andreatta} and used in the related literature. However, as suggested by a reader of these works, we could perhaps use the term {\em metagesture}.} Parametrized hypergestures are described in Section \ref{hypergesture_par}.

\subsection{Gesture morphisms}\label{gesture morphism}

A {\em morphism of two gestures} transforms one gesture into another.
Let us refer to diagram \ref{diagram_1}, containing two gestures, $g:\Delta\rightarrow \vec X$ and $h:\Sigma\rightarrow \vec Y$. We define two functions, $t$ and $m$. The function $t$ is a transformation\footnote{For example, a transformation $t^\ast$ modifies the skeleton $\Delta = .\rightarrow.$, that is an arrow between two points, into another skeleton $\Gamma = .\rightarrow.\rightarrow .$, with one more point and arrow.} from the skeleton $\Delta$ to the skeleton $\Sigma$, and the continuous function\footnote{More precisely, $m$ is a continuous functor of topological categories.} $m$, from the topological space $X$ to $Y$. We can also define $\vec m:\vec X\rightarrow \vec Y$ for the curves having points in $X$ and $Y$ respectively. These functions, jointly with the requirement of commutativity of diagram \ref{diagram_1}, define a morphism between gestures: in fact, a morphism of gestures is such a pair $f=(t,m)$ of morphisms of digraphs and topological spaces, respectively, such that $\vec m\circ g=h\circ t$ \citep{global}.
\begin{equation}\label{diagram_1}
\begin{footnotesize}
\begin{diagram}
\Delta & \rTo^g & \vec X \\
\dTo^{t} & & \dTo^{\vec m} \\
\Sigma & \rTo^h & \vec Y
\end{diagram}
\end{footnotesize}
\end{equation}
Throughout the paper, we will frequently refer to morphisms of gestures while describing the deformation of gestures with particular characteristics, such as a gesture giving a specific musical dynamic (e.g. {\em piano}) into a gesture giving a different dynamic (e.g. {\em forte}).
In general, we can define a category of gestures, with the defined morphisms as morphisms between gestures.\footnote{The category of gestures can be described as a {\em comma category} $(id,\,\vec\cdot)$
, where $id$ is the identity on directed graphs, and $\vec\cdot$
, defined from spaces to directed graphs, is the functor that associates to a space $X$ the digraph $\vec X$
with the paths in $X$.} Such a formalism can be seen in terms of $2$-categories \citep{cat1}. For example, two points in the space are 0-cells; a gesture between them is a 1-cell, and a hypergesture between two gestures is a 2-cell. If we treat point as gestures, 1-cells are hypergestures, and 2-cells are hyper-hypergestures.


\subsection{Homotopies of gestures}


In the {\em gesture follower} system, the identification of gestures is based on recognition of classes of gestures, that we can consider as homotopy classes of gestures \citep{follower, greek_2, fran, visi_miranda_2}. The authors do not use a category theoretic framework.

A homotopy of two gestures connects them via an entire family of intermediary gestures---we have a progressive `deformation' of a gesture into the other. To understand homotopies of gestures, we can intuitively think of examples such as anamorphosis frequently used in Escher's drawings. We consider a parameter $\lambda$, with $\lambda=0$ for the first gesture, and $\lambda=1$ for the second one. The progressive transformation of the first gesture into the second one is given by all the gestures labeled with values of $\lambda$ ranging between $0$ and $1$. This means that for every $\lambda\in[0,1]$, we have a gesture. A homotopic connection between gestures does not automatically imply also a functional relation between them: it is another result, not true {\em a priori}. The use of a parameter $\lambda$ for gestural comparisons recalls fuzzy logic\footnote{We will not delve into details about fuzzy logic. Whereas Definition \ref{theorem similarity} assesses a criterion of gestural similarity, we can more precisely talk about {\em degrees} of gestural similarities, with infinite intermediate values between perfect similarity of two identical gestures and two completely unrelated ones.} \citep{kosko}, as discussed in the concluding Section \ref{conclusion}.

Some more considerations about homotopy, with references to the physics of sound, are given in Paragraph \ref{hom} of the Appendix, with a ``physical'' statement about similar gestures (see the Euristic Conjecture \ref{conjecture_1}).

\subsection{Hypergestures}\label{hypergesture}

The concept of {\em hypergesture} generalizes that of homotopy.
As explained before, in a (topological) gesture space, where each point is a gesture, a curve connecting points is a {\em gesture of gestures}, a hypergesture.
The difference between hypergestures and morphisms is also stressed in \citep{mazzola_andreatta}. In the case of morphisms, we verify the relations between skeleta and bodies that are compatible; in the second case, we build curves in the space of gestures.
A particular case of hypergesture is the world-sheet, the ideal surface connecting symbolic and physical gestures, see \citep{mcm15} for a first example. However, in this article we will not use the formalism of world-sheets for our analysis.

The main structure gesture / hypergesture suggests the use of the $2$-categories formalism. In principle, because we can build hypergestures of higher order, we may extend the formalism to $\infty-$categories. The topic of infinite categories and quasicategories was investigated extensively by \cite{Joyal} and \cite{infinite}.
Music can be analyzed through nested structures. For example, a gesture (0-cell) connects two points (an arrow between points); a {\em crescendo} is a hypergesture (1-cell) that connects less loud gestures with louder gestures (an arrow between two arrows); and an {\em accelerando} is a hyper-hypergesture (2-cell) that transforms a slower crescendo into a faster crescendo (an arrow between arrows connecting arrows). We can build N-cells (and ideally $\infty$-cells) by adding each time transformations and transformations of transformations.




\section{Mathematics for similar gestures of orchestral musicians}\label{orchestra gest sim}

Let us start from a simple pianistic gesture (a primitive up-down), and an equivalent movement for a percussive gesture, let it be on a vibraphone or a xylophone. In both cases there is a `percussive' gesture/generation of the sound. 

Let $g_{\phi}$ be the physical gesture of the pianist $P$, and $h_{\phi}$ the gesture of the percussionist $P_c$. Let us for simplicity choose the same skeleton for both. The two gestures are:
\begin{equation}
g_{\phi}:\Delta\rightarrow \vec X,\,\,\,\,h_{\phi}:\Delta\rightarrow \vec Y.
\end{equation}

The topological spaces for the two instruments are different, and so is the space of their curves. We start with the same skeleton $\Delta$ for both gestures, an arrow connecting two points.
We use the formalism of category theory to characterize the action of the potential.\footnote{If restricted to the physical curve, the potential is related to the force as known in common physical situations. And, of course, the force strongly influences the touch and the final spectral result at the level of acoustics. A more general `artistic force' determines the shape of the entire surface of the world-sheet.} For a reference about world-sheet, force field, and potential see Paragraph \ref{musical forces} in the Appendix.
Let $g_{\phi}$ be the physical gesture of the pianist, represented here as the result of the action of a function $f_P:g_{\sigma}\rightarrow g_{\phi}$, where $P$ stands for the `piano' instrument. The function $f_P$ has, as its variable, symbolic gestures $g_{\sigma}$ with unspecified dynamics (that we will indicate as a $0$-potential). The physical pianistic gesture corresponding to a {\em forte} dynamic is here given as a result of the same function, with a non-vanishing potential $F_P$ as second argument: 
\begin{equation}\label{functionality}
g_{\phi}=f_P(g_{\sigma},0),\,\,\,\,g_{\phi}^{F}=f_P(g_{\sigma},F_P).
\end{equation}
Here, $g_{\phi}$ stands for $g_{\phi}^0$.
Equivalently, the percussionist's gesture is:
\begin{equation}
h_{\phi}=f_{Pc}(h_{\sigma},0),\,\,\,\,h_{\phi}^{F}=f_{Pc}(h_{\sigma},F_{Pc}).
\end{equation}
We can also represent the same situation with the (commutative) diagram \ref{diagram_3}, where $g_{\phi}^{F}=F\circ g_{\phi}$.
The operator $F$ transforms `generic' curves in $\vec X$ into {\em forte} curves in $\vec X$, that is, $(g_{\phi}\in)\vec X\rightarrow (g_{\phi}^F\in)\vec X$.
\begin{equation}\label{diagram_3}
\begin{footnotesize}
\begin{diagram}
\Delta 	&\rTo^{g_{\phi}^{F}}	&\vec X\\
	&\rdTo^{g_{\phi}} 	& \uTo_{F}\\
	&	&\vec X\\
\end{diagram}
\end{footnotesize}
\end{equation}

The digraph morphism $F:\vec X\rightarrow \vec X$ is an endomorphism of arrows of $\vec X$, because the gestural curves with unspecified dynamic (that we indicate here as with $0$-potential) are in the same space of the curves with {\em forte} dynamic ($F$-potential for us).
The digraph morphism $F:\vec X\rightarrow \vec X$ indicates a curve in the space of potentials, as already done in the case of voice \citep{ToM_III}.
We change the values, but we always are in the same space.
There is an entire family of gestures generated by the choice of potential {\em forte}. There is a homotopy of connected (linked) gestures, with the second obtained from the first via a functional composition (and this, they are not independent).

We have been composing plain morphisms of graphs. To simplify the former and the following diagrams, as well as to give a description that is coherent with the 1-cell and the 2-cell definitions, we can use 2-cell notation.
In a $2$-categorical context, diagram \ref{diagram_3} can be modified as diagram \ref{diagram_3bis}.
\begin{equation}\label{diagram_3bis}
\begin{tikzcd}
\Delta
  \arrow[r,bend left=70, "g_{\phi}"{name=U, below}] 
  \arrow[r,bend right=50, "g^F_{\phi}"{name=D}]  
& 
\vec X
\arrow[Rightarrow, from=U, to=D, "F"]
\end{tikzcd}
\end{equation} 

We are transforming the curves within the space $\vec X$, modifying the loudness coordinate via the velocity/acceleration.
We may also describe such a musical situation in the space of phases with position and speed, as often done in physics. 
In such a space, a {\em crescendo} would be easily described as a hypergesture connecting gestures with different coordinates of loudness.

To summarize, we start from a physical gesture $g_{\phi}=f(g_{\sigma},\circ)$ from a certain symbolic gesture $g_{\sigma}$ and a potential $\circ$, giving the dynamic. We could think of $\circ$ as a generic potential $V$, that is then specialized when a particular dynamic ({\em forte}, {\em mezzoforte}) is indicated.
A transformation $0\rightarrow F$ with the {\em forte} induces a deformation of the gestural curve into another realizing the {\em forte} dynamic as result of the interaction between performer and musical instrument: we have $g_{\phi}\rightarrow g_{\phi}^F$, where $g_{\phi}^F=f(g_{\sigma},F)$. The transformation of the potential $0\rightarrow F$ can be described by a continuous function\footnote{More generally, we can define a gesture via a function $f$ containing such a $V(t)$ function.
\[
f(g_{\sigma},V(t))=:g_{\phi}(t)= \left\{
\begin{array}{lr}
g_{\phi}^{0}, & for\, t=0\\
g_{\phi}^{F}, & for\, t=1\\
\end{array}\right\}
\]}
 $V(t)$, that is equal to $0$ at $t=0$,
 and to $F$ at $t=1$, the coordinate $t$ being a generic parameter, or a time coordinate.
We can describe a crescendo via diagram \ref{diagram_4}. We start from a gesture $g$ with unspecified dynamic, we then identify the required transformation (i.e. gestural deformation) to get the {\em piano} and {\em forte} respectively, and {\em then} we define a {\em temporal} transformation $PF$ that brings the {\em piano} gesture $g_{\phi}^{P}$ into the {\em forte} gesture $g_{\phi}^{F}$, that {\em is} the crescendo. In terms of potentials, we have the transformation $F={PF}\circ P$.
\begin{equation}\label{diagram_4}
\begin{footnotesize}
\begin{diagram}
g_{\phi}^{P} &	&\rTo^{PF} &&	g_{\phi}^{F} \\
&\luTo^{P}	&	&\ruTo^{F}	&\\
&&		g_{\phi}	&&\\
\end{diagram}
\end{footnotesize}
\end{equation}
In terms of $2$-category formalism, diagram \ref{diagram_4} can be re-drawn as diagram \ref{diagram_4bis}, using the vertical composition property.
\begin{equation}\label{diagram_4bis}
\begin{tikzcd}
\Delta
  \arrow[r,bend left=80, "g_{\phi}"{name=U, above}] 
    \arrow[r, "g^P_{\phi}"{name=N}] 
  \arrow[r, bend right=80, "g^F_{\phi}"{name=D, below}]  
& 
\vec X
\arrow[Rightarrow, from=U, to=N, "P"]
\arrow[Rightarrow, from=N, to=D, "PF"]
&=&
\Delta
  \arrow[r,bend left=70, "g_{\phi}"{name=U, above}] 
  \arrow[r,bend right=50, "g^F_{\phi}"{name=D, below}]  
& 
\vec X
\arrow[Rightarrow, from=U, to=D, "F"]
\end{tikzcd}
\end{equation}


Let us now explicitly use functionality.
We reach the {\em forte} curves in $\vec X$ via $f\left(\circ, F\right)$ for the {\em forte} $F$, and $f(\circ,0)$ without any specified dynamic, where the symbol $\circ$ indicates the place reserved for the symbolic gesture.
Here, we choose an elementary skeleton\footnote{We denote by $\uparrow$ the digraph having two vertices that are connected by one arrow.} $\Delta = \uparrow$, and a simple symbolic gesture that is the image of such an abstract curve, indicated here also with $\uparrow$.
We have  $f\left(\uparrow, F\right)=g_{\phi}^{F}$ and $f(\uparrow,0)=g_{\phi}$, $F\circ f(\uparrow,0)=f(\uparrow, F)$, and the diagram \ref{diagram_5} commutes.
\begin{equation}\label{diagram_5}
\begin{footnotesize}
\begin{diagram}
\uparrow	&\rTo^{f(\uparrow, F)}	&\vec X\\
		&\rdTo^{f(\uparrow,0)} 	& \uTo_{F}\\	
		&					&\vec X\\
\end{diagram}
\end{footnotesize}
\end{equation}
Two main questions arise:
\begin{enumerate}
\item Is this homotopy?
\item Is the process $g_{\phi}\rightarrow g_{\phi}^{F}$ functional (well-defined)?
\end{enumerate}
The answer to both questions is affirmative for the following reasons.
\begin{enumerate}
\item Homotopy. There is homotopy between the two gestures, because there is an entire family of intermediate gestures continuously connecting the first and the second one. The potentials themselves are connected by a family of potentials, going from $0$ and reaching $F$.

\item Functionality. The process is also functional, for the reason expressed above: the gesture $f(\uparrow, F)$ is functionally obtained from $f(\uparrow,0)$ via the composition with $F$. In fact, $g_{\phi}^{F}\equiv F\circ g_{\phi}$.
\end{enumerate}

Let us now extend the previous discussion to include percussionists' gestures. We again use the label $P$ for the piano, and the label $Pc$ for the percussion. In fact, this time we distinguish between $F_P$ and $F_{Pc}$, acting on different spaces.
Diagram \ref{diagram_6} shows the pianist's gesture (right) and the percussionist's gesture (left), having the same skeleton $\Delta$.
\begin{equation}\label{diagram_6}
\begin{footnotesize}
\begin{diagram}
\vec Y &\lTo^{h_{\phi}^F} & \Delta 	&\rTo^{g_{\phi}^F}	&\vec X \\
\uTo^{F_{Pc}}  	& \ldTo^{h_{\phi}} & &\rdTo^{g_{\phi}} 	& \uTo_{F_P}\\
\vec Y	& & &	&\vec X\\
\end{diagram}
\end{footnotesize}
\end{equation}
\begin{equation}\label{diagram_6bis}
\begin{tikzcd}
\vec Y
&
\Delta
    \arrow[r,bend right=60, "g^F_{\phi}"{name=D}]  
    \arrow[r,bend left=70, "g_{\phi}"{name=U, above}] 
  \arrow[l,bend left=60, "h^F_{\phi}"{name=Z}]
   \arrow[l,bend right=70, "h_{\phi}"{name=V, above}]  
   \arrow[Rightarrow, from=V, to=Z, "F_{Pc}"]
&
  \arrow[Rightarrow, from=U, to=D, "F_P"] 
\vec X
\end{tikzcd}
\end{equation}
Diagram \ref{diagram_6bis} shows the same content of diagram \ref{diagram_6} in a $2$-categorical context.\footnote{We can observe that, in the framework of $2$-category formalism, vertical and horizontal composition properties are intuitively verified. Vertical composition means, in our orchestral context, transformation of loudness: a {\em piano} gesture can be deformed into a {\em forte} gesture, that can be deformed into a {\em fortissimo} gesture. Horizontal composition, given the same skeleton, means here transition from the space of gestures for a musical instrument, for example piano, to the space of gestures of another instrument, as percussion, and then to another instrument, such as violin.}

The two triangles in diagram \ref{diagram_6} show functional transformations. How can we connect $\vec Y$ to $\vec X$ and compare the two gestures $g$ and $h$? We think of $g$ and $h$ as being `gesturally similar,' because we can transform the one into the other via an homotopic transformation, {\em and} they lead to similar timbres both modified into {\em forte}---see the Euristic Conjecture \ref{conjecture_1}. We will now consider two different skeleta, and the change of skeleta $\Delta\rightarrow\Gamma$, as shown in diagram \ref{diagram_7}, provided we are given a continuous function $m$ such that $m:X\rightarrow Y$.
\begin{equation}\label{diagram_7}
\begin{footnotesize}
\begin{diagram}
&	&		&	&	&\vec X\\
&	&		& 	&\ruTo(4,2)^{g_{\phi}}\ldTo^{F_P}&\dTo_{\vec m}\\
&\Delta	&\rTo^{g_{\phi}^F}		&\vec X	&\\
&\dTo{t}	&		&\dTo{\vec m} 		&\\
&\Gamma	&\rTo^{h_{\phi}^F}		&\vec Y	&\\
&	&\rdTo(4,2)^{h_{\phi}}	& 	&\luTo^{F_{Pc}}&\\
&	&		&		& &\vec Y
\end{diagram}
\end{footnotesize}
\end{equation}
The two triangles are commutative because $g_{\phi}^F=F_P\circ g_{\phi}$ and $h_{\phi}^F=F_{Pc}\circ h_{\phi}$, as well the square of diagram \ref{diagram_8}.
\begin{equation}\label{diagram_8}
\begin{footnotesize}
\begin{diagram}
&\Delta	&\rTo^{g_{\phi}^F}		&\vec X	&\\
&\dTo{t}	&		&\dTo{\vec m} 		&\\
&\Gamma	&\rTo^{h_{\phi}^F}		&\vec Y	&\\
\end{diagram}
\end{footnotesize}
\end{equation}
This means that there is a morphism of gestures, as defined in Section \ref{gesture morphism}.
Diagram \ref{diagram_7}, with $2$-category formalism, looks like diagram \ref{diagram_7bis}.

\begin{small}
\begin{equation}\label{diagram_7bis}
\begin{tikzcd}[swap,bend angle=45]
\Delta
\arrow{ddd}{t}
    \arrow[r,bend right=70, "g_{\phi}^F"{name=D}]  
    \arrow[r,bend left=80, "g_{\phi}"{name=U, above}] 
&
  \arrow[Rightarrow, from=U, to=D, "F_P"] 
\vec X
\arrow{ddd}{\vec m}
   \\
   \\
   \\
   \Gamma
  \arrow[r,bend right=70, "h_{\phi}"{name=V}]
   \arrow[r,bend left=80, "h_{\phi}^F"{name=Z, above}]  
&
\vec Y
   \arrow[Rightarrow, from=V, to=Z, "F_{P_c}"]
\end{tikzcd}
\end{equation}
\end{small}



Under which conditions is diagram \ref{diagram_7} commutative also with $\vec X$ and $\vec Y$? 
The similarity definition also includes the commutativity, in diagram \ref{diagram_7}, of the square shown in diagram \ref{diagram_new}.
\begin{equation}\label{diagram_new}
\begin{footnotesize}
\begin{diagram}
&\vec X	&\lTo^{F_P}		&\vec X	&\\
&\dTo{\vec m}	&		&\dTo{\vec m} 		&\\
&\vec Y	&\lTo^{F_{Pc}}		&\vec Y	&\\
\end{diagram}
\end{footnotesize}
\end{equation}

Similar gestural generators are supposed to provoke similar changes in the timbres, for example the creation of two {\em forte} gestures from gestures with an unspecified loudness. But, first of all, the changes have to happen in the appropriate spaces.
Diagram \ref{diagram_7} is not commutative when, for example, the generative process is similar but the resulting spectrum is not. If we strongly hit the flute with a piano hammer, it will not sound louder, and a strong use of the bow on a flute will not make a stronger flute sound. They are examples of {\em forte} gestures in non-appropriate spaces for the flute. This means that we cannot apply the same ``generator'' $F_P$ of pianistic {\em forte}, to obtain a {\em forte} gesture on flute. We need another generator, a $F_{Fl}$, to transform generic curves of flute playing into {\em forte} curves (gestures having as result a {\em forte} sound). $F_{Fl}$ does not act on $\vec X$ but on some $\vec Z$, where $\vec Z$ is the space of curves for the flute playing. The spaces $\vec X$ and $\vec Z$ can be connected by a suitable $m$ function to transform the space of gestures of an instrument into the space of gestures of another instrument.

We can finally give a first definition of gestural similarity, as expressed in Definition \ref{theorem similarity}. 
\begin{definition}\label{theorem similarity}
Two gestures $g_{\phi}^F:\Delta\rightarrow\vec X$ and $h_{\phi}^F:\Gamma\rightarrow\vec Y$ are {\em similar}, and we write $g_{\phi}^F\sim h_{\phi}^F$, if it is possible to find:
\begin{enumerate}
\item A morphism from the first to the second, it means, two functions $t:\Delta\rightarrow\Gamma$ and $m:X\rightarrow Y$ such that diagram \ref{diagram_8} is commutative;
\item Two homotopic transformations\footnote{Two {\em homotopic} curves can be continuously transformed the one into the other. A homotopic transformation continuously transforms curves in curves.} $F_P:\vec X\rightarrow\vec X$ and $F_{Pc}:\vec Y\rightarrow\vec Y$, with the {\em common generators} $F_P$ and $F_{Pc}$ {\em acting in the same way}\footnote{Similar deformations in their respective spaces lead to similar effects in their resulting sound spectra. See Paragraph \ref{hom} in the Appendix for two graphic representations, Figures \ref{non_sim} and \ref{3d_diagram}.} in their spaces respectively,
such that $F_P\circ g_{\phi}=g_{\phi}^F$ and $F_{Pc}\circ h_{\phi}=h_{\phi}^{F}$.
\end{enumerate}
\end{definition}
Let us see more details in Example \ref{example_sim}.
\begin{example}\label{example_sim}

Let us suppose {\em ad absurdum} that two gestures $g_1,\,g_2$ for the same instrument (topological space $\vec X$), with unspecified dynamic but sharing the same skeleton, are similar. If the skeleta are the same, we have an identity map between them. We also have an identity for the topological space $X\rightarrow X$, being the two gestures  for the same instrument. So we have a gesture morphism with identities in this case, and point $(1)$ is satisfied. However, because the dynamic is unspecified, we also have unspecified generators. This implies that point $(2)$ is not satisfied, and thus the hypothesis of similarity is not verified.
In fact, two gestures cannot be similar if they are obtained via two generators acting in different ways in their respective spaces (the same in our case).
We can consider a {\em forte} gesture and a {\em piano} gesture, both for the same musical instrument, for example percussion $Pc$, with the same skeleton. We have an identity for the topological space $Y\rightarrow Y$. The two generators $F_{Pc}:\vec Y\rightarrow\vec Y$ and $P_{Pc}:\vec Y\rightarrow\vec Y$ do not act in the same way: the first selects {\em forte} curves, while the second {\em piano} curves. We have $F_{Pc}\circ h_{\phi}=h_{\phi}^F$ and $P_{Pc}\circ h_{\phi}=h_{\phi}^P$, and point $(2)$ is not satisfied.
\end{example}
See also the Heuristic Conjecture \ref{conjecture_1} in Paragraph \ref{hom} of the Appendix for a physical condition about acoustical spectra. The musical meaning of Definition~\ref{theorem similarity} is the following.
A {\em forte} gesture for the pianist and a {\em forte} gesture of the percussionist (in our example) are {\em similar} because:
\begin{enumerate}
\item It is possible to transfer (connect) the {\em forte} on the piano to the {\em forte} on a percussion; 
\item The generator that transforms generic (neutral) gestures on piano into {\em forte} gestures on piano is the same generator that transforms generic gestures on a percussion into {\em forte} gestures on a percussion, acting in different spaces in the same way.
\end{enumerate}
Both {\em forte} gestures\footnote{We discussed loudness, but gestural similarity can also involve articulation, and in some cases also rhythmic-melodic profiles. Even specific harmonic sequences can suggest particular gestural solutions. We can think of a deceptive cadence, highlighted by the performer with a {\em fermata} or a {\em forte}. In fact, elements from musical analysis can act as weights for gestures.
} can be projected into the {\em forte} conducting gestures,
using the categorical notion of {\em colimit} (Section \ref{role conductor}).
To decide if the two generators $F_P$ and $F_{Pc}$ act in the same way in their spaces, and, consequently, if two resulting gestures are similar, we can also use concepts from fuzzy logic. A gesture can belong to a category but also partially to another one. For example, a range of variability allows the distinction of several {\em forte} gestures, that are all equally well-working for the same musical passage. This opens another potential research in the application of fuzzy logic to category theory \citep{fuzzy1} for music.

\subsection{Conducting gestures}\label{role conductor}

The conductor's gesture is specialized into gestures of all orchestral musicians, and it may be described as their {\em colimit}. For the listener, the conducting gesture can constitute a {\em limit}.\footnote{Intuitively, {\em limits} and {\em colimits} are generalizations of {\em products} and {\em coproducts}, respectively. The product is a special case of the limit, with discrete categories. The coproduct (also called {\em sum}) is the dual of the product, obtained by reversing the arrows. Given an object $P$ and two maps $p_1:P\rightarrow B_1,\,p_2:P\rightarrow B_2$, $P$ is a product of $B_1,\,B_2$ if for each object $X$ and for each pair of arrows $f_1,\,f_2$ we have only and only one arrow $f:X\rightarrow P$ such that $f_1=p_1f,\,f_2=p_2f$ \citep{lawvere}, see diagram \ref{product}.
\begin{equation}\label{product}
\begin{footnotesize}
\begin{diagram}
&	&&		\mbox{X }		& &\\
&	&\ldTo(2,4)^{f_2}&\dTo	 ^{f}&			\rdTo(2,4)^{f_1}	& &\\
&	&&			P	& &\\
&			& \ldTo^{p_2}&&\rdTo^{p_1}		 &  \\
&			B_2&			&	& &B_1 &    \\
\end{diagram}
\end{footnotesize}
\end{equation}}
We can think of a person trying to get information about the orchestral sound by watching conducting gestures on a video with the audio turned off.
We are considering limit and colimit because the listener (i. e., each of the listeners in the audience) refers to the orchestra, and all the orchestral instruments refer to the conductor. Moreover, the conductor, while studying the score, prepares his or her gesture depending on the orchestration, meaning what orchestral musicians are supposed to play. The listener creates his or her mental idea of the music depending on what orchestral performers are playing.

We indicate orchestral gestures as $\mathcal{D}$, conductor's gestures is their colim$(\mathcal{D})$, while the audience is the lim$(\mathcal{D})$.

We can draw diagram \ref{diagram_9}, including the contribution of the conductor (with $\vec C$), where
$F_C$ ($\vec m_{C,F}$), associated to the conductor's gesture, transforms conducting generic curves (with unspecified dynamic) into {\em forte} curves.
In this diagram, to simplify the graphical representation, we omit the skeleton of the conducting gesture.
\begin{equation}\label{diagram_9}
\begin{footnotesize}
\begin{diagram}
&	&		&	&	&\vec X	&	&	\\
&	&		& 	&\ruTo(4,2)^{g_{\phi}}\ldTo^{F_P}&	&\rdTo^{\vec m_{C,P,0}}&		\\
&\Delta	&\rTo^{g_{\phi}^F}		&\vec X	&	&& &\vec C		\\
&\dTo{t}	&		&\dTo{\vec m} 		&\rdTo^{\vec m_{C,P,F}}	
&&\ldTo^{F_C}	&\dTo^{Id}	\\
&	&		& 		&		&\vec C	\\
&	&		& 		&\ruTo^{\vec m_{C,Pc,F}}		&&\luTo^{F_C}	\\
&\Gamma	&\rTo^{h_{\phi}^F}		&\vec Y	&		&&&\vec C	\\
&	&\rdTo(4,2)^{h_{\phi}}	& 	&\luTo^{F_{Pc}}&	&\ruTo^{\vec m_{C,Pc,0}}		\\
&	&		&		& &\vec Y	&	&		
\end{diagram}
\end{footnotesize}
\end{equation}
In diagram \ref{diagram_9}, we distinguish the functions bringing the conductor's gesture to the pianist and to the percussionist. Summarizing, we have:
\begin{itemize}
\item gesture without specified dynamic of the pianist reflected into the basic metric gesture (no indication of dynamic) of the conductor: $\vec m_{C,P,0}: \vec X\rightarrow \vec C$,
\item gesture without specified dynamic of the percussionist reflected into the basic metric gesture (no indication of dynamic) of the conductor: $	\vec m_{C,Pc,0}:\vec Y\rightarrow \vec C$,
\item gesture of {\em forte} for the pianist reflected into the conductor's {\em forte}: $\vec m_{C,P,F}: \vec X\rightarrow \vec C$,
\item gesture of {\em forte} for the percussionist reflected into the conductor's {\em forte}: $\vec m_{C,Pc,F}: \vec Y\rightarrow \vec C$,
\item basic metric gesture transformed into the {\em forte} one of the conductor: $F_{C}:\vec C\rightarrow \vec C$,
\item identity $Id$ transforming conductor's basic metric gesture into itself. 
\end{itemize}
Diagrams \ref{diagram_10} show the two commutative squares conductor-pianist and conductor-percussionist from diagram \ref{diagram_9}, respectively.
\begin{equation}\label{diagram_10}
\begin{footnotesize}
\begin{diagram}
\vec X & \rTo^{\vec m_{C,P,0}} & \vec C \\
\dTo^{F_P} & & \dTo^{F_C} \\
\vec X & \rTo^{\vec m_{C,P,F}} & \vec C
\end{diagram}
\end{footnotesize}
\begin{footnotesize}
\hspace{40pt}
\begin{diagram}
\vec Y & \rTo^{\vec m_{C,Pc,0}} & \vec C \\
\dTo^{F_{Pc}} & & \dTo^{F_C} \\
\vec Y & \rTo^{\vec m_{C,Pc,F}} & \vec C
\end{diagram}
\end{footnotesize}
\end{equation}
These diagrams commute because it is equivalently possible that:
\begin{itemize}
\item The pianist plays {\em forte} because the conducting basic metric gesture has been modified by the operator {\em forte};
\item The conductor gives a basic metric gesture and the pianist decides to play {\em forte} for reasons of expressivity or written indication.
\end{itemize}
This second case can happen, for example, for scores where rhythm and other specific
indications are so complicated that the conductor needs to give simple metric indications. In formulas, we have, for the pianist, $F_C\circ\vec m_{C,P,0}= \vec m_{C,P,F}\circ F_P$, and the same for the percussionist, substituting $P_c$ to $P$.

Diagram \ref{diagram_12} shows the gestural similarity between the {\em forte} gesture of the pianist and the percussionist (via $\vec m$), and their reflection into the conducting {\em forte} gesture, for the pianist (via $\vec m_{C,P,F}$) and the percussionist (via $\vec m_{C,Pc,F}$).
\begin{equation}\label{diagram_12}
\begin{footnotesize}
\begin{diagram}
\vec X &  \\
\dTo^{\vec m}  & \rdTo^{\vec m_{C,P,F}} \\
&  &\vec C \\
&   \ruTo^{\vec m_{C,Pc,F}} \\
\vec Y
\end{diagram}
\end{footnotesize}
\end{equation}

Using $2$-categories, we can schematize diagram \ref{diagram_9} as shown by diagram \ref{diagram_9bis}, with functors transforming $\vec m_{C,P,0}$ in $\vec m_{C,P,F}$ and $\vec m_{C,P_c,0}$ in $\vec m_{C,P_c,F}$, respectively.\footnote{This is not shown in the diagram.}

\begin{small}
\begin{equation}\label{diagram_9bis}
\begin{tikzcd}[swap,bend angle=45]
\Delta
\arrow{ddddd}{t}
    \arrow[r,bend right=70, "g_{\phi}^F"{name=D}]  
    \arrow[r,bend left=80, "g_{\phi}"{name=U, above}] 
&
  \arrow[Rightarrow, from=U, to=D, "F_P"] 
\vec X
\arrow{ddddd}{\vec m}
\arrow[dddrrr,bend right=30, "\vec m_{C,P,0}"{name=A}]
\arrow[dddrrr,bend left=30, "\,\,\,\,\,\,\,\,\,\,\,\,\vec m_{C,P,F}"{name=B, above}]
   \\
   \\
   \\ & & & & \vec C
   \\ 
   \\ 
   \Gamma
  \arrow[r,bend right=70, "h_{\phi}^F"{name=Z}]
   \arrow[r,bend left=80, "h_{\phi}"{name=V, above}]  
&
\vec Y
\arrow[uurrr,bend right=30, "\vec m_{C,P_c,0}"{name=C}]
\arrow[uurrr,bend left=30, "\vec m_{C,P_c,F}",{name=H, above}]
   \arrow[Rightarrow, from=V, to=Z, "F_{P_c}"]
\end{tikzcd}
\end{equation}
\end{small}

We can apply such a diagram to all musicians of the orchestra, getting the {\em colimit}: orchestra $= \mathcal{D}\longrightarrow$ conductor $= \mbox{colimit} (\mathcal{D}$).
We can define the conductor-orchestra-listener via diagram \ref{diagram_13} in terms of limits and colimits.

\begin{equation}\label{diagram_13}
\begin{footnotesize}
\begin{diagram}
& &conductor=\mbox{colim}(\mathcal{D}) \\
&	&\uTo\uTo\uTo	&	\\
&	&orchestra=\mathcal{D}	&	\\
&	&\uTo\uTo\uTo	&	\\
& &listener=\lim(\mathcal{D})
\end{diagram}
\end{footnotesize}
\end{equation}
If we schematically represent the orchestral gestures as in equation \ref{diagram_14},
we can define the limit as the listener,
and the colimit as the conductor,
see diagram \ref{diagram_15}. The choice of listener as a limit, even if considered as a metaphor, satisfies the universal property because all the ``listening and perception activities'' can be reduced to the listener, in contraposition to the sound-production activities. The conductor plays the opposite role: all the ``sound production'' gestural activities can be related to the conducting gesture, that is a pure gesture without any direct sound production.
\begin{equation}\label{diagram_14}
\mathcal{D}=
\begin{diagram}
&			X_{\lambda}&			&\rTo^{f_{\lambda\mu}}	& &X_{\mu} &    \\
&	&\luTo^{f_{\kappa\lambda}}		&	&\ruTo^{f_{\kappa\mu}}	&	\\
&		&			&X_{\kappa}	& & &    \\
\end{diagram}
\end{equation}

\begin{equation}\label{diagram_15}
\begin{footnotesize}
\begin{diagram}
&	&&		\lim\mathcal{D}		& &\\
&	&\ldTo(2,4)&\uTo	 ^{!}&			\rdTo(2,4)	& &\\
&	&&			Z	& &\\
&			& \ldTo^{\pi_{\lambda}}&&\rdTo^{\pi_{\mu}}		 &  \\
&			X_{\lambda}&			&\rTo^{f_{\lambda\mu}}	& &X_{\mu} &    \\
&	&\luTo^{f_{\kappa\lambda}}		&	&\ruTo^{f_{\kappa\mu}}	&	\\
&		&			&X_{\kappa}	& & &    \\
\end{diagram}
\hspace{40pt}
\begin{diagram}
&	&&		\mbox{colim }\mathcal{D}		& &\\
&	&\ruTo(2,4)&\dTo	 ^{!}&			\luTo(2,4)	& &\\
&	&&			Z	& &\\
&			& \ruTo^{\pi'_{\lambda}}&&\luTo^{\pi'_{\mu}}		 &  \\
&			X_{\lambda}&			&\rTo^{f_{\lambda\mu}}	& &X_{\mu} &    \\
&	&\luTo^{f_{\kappa\lambda}}		&	&\ruTo^{f_{\kappa\mu}}	&	\\
&		&			&X_{\kappa}	& & &    \\
\end{diagram}
\end{footnotesize}
\end{equation}


We can also specify the colimit diagram as shown by diagram \ref{diagram_16}, where $\mathcal{D}$ is the orchestra, and the colimit {\em is} the conductor's gesture $\Delta_{conductor}\rightarrow\vec X_{conductor}$.
\begin{footnotesize}
\begin{equation}\label{diagram_16}
\begin{diagram}
& & & & & & & \Gamma & \rTo^g & \vec Z \\
& & & & & & & \ruTo(3,2)^{\exists!} & & \luTo^{\forall} \\
& & & & \Delta_{conductor} &\rTo_{g_{conductor}} &  \vec X_{conductor} & \rTo& & &\mathcal{D}
\end{diagram}
\end{equation}
\end{footnotesize}

Let us now include the detail of (some) orchestral gestures.
In diagram \ref{diagram_17}, for simplicity we represent just one arrow from one orchestral gesture, $g_\lambda$, to the conductor; however, {\em each} orchestral gesture, it means each $g_i$ ($g_{\lambda}$, $g_{\kappa}$, $g_{\mu}$) in the diagram \ref{diagram_17}, has an arrow to the conductor's gesture $g_{conductor}$, and also to the $g$ gesture.

\begin{equation}\label{diagram_17}
\begin{footnotesize}
\begin{diagram}
& & & & & & & \Gamma & \rTo^g & \vec Z \\
& & & & & & & \ruTo(3,2)^{\exists !} &  \uTo^{\forall} \\
& & & & \Delta_{conductor} &\rTo_{g_{conductor}} &  \vec X_{conductor} & & & && & \\
&&&& & & &\Delta_{\kappa}&  \rTo^{g_{\kappa}} & \vec X_{\kappa} & \\
&&& & &\uTo&   \ldTo(3,3) & &\ldTo(3,3) & \\
 &&& & & &   &\dTo & &\dTo & \\ 
 &&&& \Delta_{\lambda} &  \rTo^{g_{\lambda}} & \vec X_{\lambda}\\
&&&& &  \rdTo(3,3) &  & \rdTo(3,3) & \\
&&&& & & & \dTo & &    & \\
&&&& & & &\Delta_{\mu} & \rTo^{g_{\mu}} & \vec X_{\mu} \\
\end{diagram}
\end{footnotesize}
\end{equation}
We will not deal here with a description of the universal properties of these limits and colimits, preferring some more comments on their musical meaning and implications. The gestures of the orchestral musicians can be (mathematically) injected into the conductor's ones.

Intuitively, we would say the opposite: the gestures of the conductor are the initial source of movement for the orchestra, and the listener is the final target who sees the gestures, and listens to their sound result.
This means that the conductor's movements are specified into the gestures of orchestral musicians, and collected as a whole result by the listener. This would imply an inverse order of the arrows in diagrams \ref{diagram_14} and \ref{diagram_15}.

However, we can also say that the listener (emotionally) projects his or her thought inside the sound and the gestures, so the arrows come from the listener and reach the orchestral gestures. We may see the orchestral gestures as contained ``inside'' the conductor's ones. 
Let us think of the violinist's, or pianist's gestures. They are much more complicated than the conducting ones. The complicated gesture of each performer refers to, and injects into the simpler gesture, i.e. the conducting one. This means that it projects into its {\em simplification}. We can interpret piano or violin gestures as something that can be envisaged into the conducting gestures.
Orchestral gestures are the development of the conducting gestures.
The listener may imagine the orchestral gestures---and the shape of music---by observing the conducting gestures, even without listening to the sound.
In fact, the conductor's gesture is the {\em terminal} element: there is one and only one function going from it to the $g:\Gamma\rightarrow \vec Z$ having the property that {\em all} orchestral gestures are injecting into it. There is a morphism of each instrument to the conductor, and there is one and only one morphism of gestures of the conductor in $Z$ that makes the diagram commutative. We can talk metaphorically about categorical adjunction between these two perspectives. A mediation between the two poles would involve the definition of exchanges and mutual feedback between performers, conductor, and listeners; the conductor is also a listener and the listener can influence the performance through applause or boos. An intermediate figure between who {\em gives indication to make} and who {\em listens to} the sounds may be a conductor that changes the gesture according to the performers' outcome, or, better, a performer of an electronic instrument that influences the sound production through the gesture itself (as a theremin player), and adjusts the gesture according to the sound.

If we opt for the second description or for the intuitive one, both listener and conductor are still at the opposite side, because:
\begin{itemize}
\item there is only one element (initial or terminal, depending on the chosen description) that collects all the orchestral gestures (for example, a ``simplified version'' of orchestral gestures with collected analogous movements), and there is one and only one morphism to the conductor's {\em motor} gesture;
\item there is only one element (again, terminal or initial, depending on the chosen description) that reaches all the orchestral gestures (for example, a more prominent sound result given by the most relevant musical structures, elements, phrases, timbres), and there is one and only one morphism to the listener's {\em perceptive} gesture.
\end{itemize}
Both listener and conductor are silent. The conductor is making active gestures, the listener is interpreting the music heard via making similarities with his or her background, past listening experiences, musical knowledge, and personal sensibility. 
The conductor belongs to the {\em motor} world, while the listener to the {\em perceptive} world. If, through a software, the listener can manipulate an electronic conductor, he or she is also a conductor, and if the conductor changes his or her gestures depending on the {\em listened sound}, he or she is also a listener. In each case, listener and conductor can only be compared with initial and terminal objects in the world of the orchestral performance.\footnote{This happens if we already have a musical score. If, in such a description, we include the composer of the orchestral score, then the colimit role may be envisaged into the composing activity. We may say that, in this extended description, the initial and final points are both in the mind: the mind of the composer, and the mind of the listener.} In the case of a piano improviser who is listening to his or her own music, listener's, (unique) performer's, and composer's roles all coincide.

Within the formalism of $2$-category,
we can generalize the description of listener $\rightarrow$ conductor $\rightarrow$ orchestra in terms of $2$-colimits and $2$-limits. In this case, horizontal composition properties (among orchestral instruments, among different conductors, among different listeners) and vertical composition properties (among gestures leading to different dynamics, articulations, tempo changes) are still valid. This can be the topic of further research in itself.





\section{Gestures from the composer to the conductor}\label{adjoint}

We can schematize the gestures of the orchestral performance as a category, and the gestures hidden in the score as another category. We can define a {\em Performance Functor} connecting symbolic to physical gestures,\footnote{The connection of single curves from (symbolic) systems of continuous curves to single curves from other (physical) systems of continuous curves, their skeleta being the same, has been investigated via branched graphs and branched world-sheets \citep{mcm15, global}.} and a {\em Compositional Functor}, for the inverse movement, from the improvisation to the music transcription/composition.
In diagram \ref{diagram_18}, we consider different skeleta, using labels $\sigma,\,\phi$ to denote symbolic and physical gestures, respectively.
Let us use the concept of functor as morphism between categories. The Performance Functor (P. F.) connects symbolic gestures to physical gestures, for piano and percussion.
\begin{equation}\label{diagram_18}
\begin{footnotesize}
\begin{diagram}
\Delta_{\sigma}& \rTo & \vec X_{\sigma} \\
\dTo & &\rdTo(3,4)^{P.F._{P}} \dTo \\
\Gamma_{\sigma}& \rTo &\vec Y_{\sigma}\\
& & \rdTo(3,4)^{P.F._{Pc}} & & \\
& & &\Delta_{\phi}& \rTo & \vec X_{\phi} \\
& & &\dTo & & \dTo & & \\
& & &\Gamma_{\phi}& \rTo &\vec Y_{\phi}\\
\end{diagram}
\end{footnotesize}
\end{equation}

\section{Gestural similarity between music and visual arts}\label{visual}

The concept of gestural similarity can be used not only in music, but also between music and visual arts.\footnote{Category theory includes functors: so, we could provocatively talk about {\em functorial aesthetics}.}  In art history, there are several examples where visual artists got inspiration from music and vice versa: we can think of Boris Mussorgsky, Luigi Russolo, and Morton Feldman, as few names.\footnote{In future developments of gestural similarity analysis, we may try to find analogies between a painting and the music inspired from it. This may be part of a more general approach to artistic movements, finding the connections between music, visual art, and poetry within a specific movement, in terms of basic ``shared gestures.''}  In movies, particular scenes are often emphasized by musical gestures: we can think of the correspondence between the frightening bow's strikes and the knife scene in {\em Psycho}.

Several studies highlight the correspondence between music and movement \citep{zbikowski}, especially with the applications in the field of soundtracks. In general, we can say that ``music is a gesture,'' it is the ``result of forces, impulses, experiences... both from individual and general'' \citep{sessions2}.

Finally, we can discuss examples of sonification (mapping of non-sound data into sound-data), and we can argue that they are more effective when there are gestural similarities: a rising shape may be described via a rising pitch sequence \citep{roffler}, or a ``staccato'' musical sequence may be compared with a collection of points on canvas. It is sufficient to see a drawing, for example, as the result of a drawing gesture. Studies in the field of psychology \citep{zbikowski}, crossmodal correspondences \citep{spence, ohala, berlin}, iconicity in linguistics \citep{mayer}, and audio-visual objects and the theory of indispensable attributes \citep{kubovy} can support these ideas. Finally, we can re-read examples of sonification \citep{music_image_book} in light of gesture theory, and we can compose new music \citep{phd_mannone}.

We can use these ideas to investigate the connections between vocal gestures in speaking and singing, and visual shapes, framing within category theory classic perceptual experiments \citep{uznadze, nobile, nobile2}.

Further developments of such a research will also include perceptual experiments to test the expected degree of similarity between music and visuals.

\section{Conclusion}\label{conclusion}

In this paper we contextualized in a categorical framework the analogies between gestures of different musicians, the gestural communication between conductor and performers. Such a formalism can also be used to compare gestures in music and in visual arts.

We proposed a first mathematical definition of gestural similarity, to allow comparisons between gestures belonging to different spaces (such as piano's and percussion's spaces) having
similar characteristics, for example all realizing a {\em forte} sonority.
Further research can address the more general question of the gesture classification.
We can use concepts from {\em fuzzy logic}, to quantify the {\em degree of similarity} between two different gestures. Moreover, fuzzy logic can intervene in a more fundamental way, to connect symbolic and physical gestures. In fact, analogously to what we do when we think of ``the'' circle (as a Platonic idea) without being able to draw it, except in terms of a fuzzy circle `not perfectly circular,' we can write a score and not be able to play {\em exactly what it is written}. Each performance can be described in fuzzy terms: $0$ for the score, $1$ for the ``perfect'' and complete musical performance, and values within $[0,1]$ to denote the performance in progress. Transition from symbolic to physical gestures can be reformulated in terms of progressive action of a fuzzy function.

About category theory, further research may explore in detail the universal properties of $2$-limits and $2$-colimits between conductor, orchestra and listener, as well as theoretical extensions and artistic implications in the field of $\infty$-categories. 

About the relations between gestural similarity and sound, future research can also involve collaboration with scholars in the field of neurosciences, more precisely about interactions between symbolic gestures and words \citep{gent1}, spoken language and arm gestures \citep{gent2}, gesture, sound, imitation in phoneme production \citep{gent3, gent4, gent5}.

In conclusion, future studies can deepen and consolidate the connections between physics, psychophysics, visual art, linguistics, and neuroscience, via a common mathematical model.

\section{Appendix: More details on the mathematical theory of musical gestures}

\subsection{Homotopy and gestures}\label{hom}

We will clarify the concept of homotopy and gestures used in this article.
Let us suppose we have two gestures, one for the piano, and the other for the violin, with the same skeleton, see Figure \ref{non_sim}. A {\em forte} gesture for the pianist is similar to the {\em forte} for the violinist. The {\em forte} is transformed into {\em piano} via a {\em diminuendo}, represented by the action of the functor $FP$ in the $2$-category (more precisely, we have $FP_P$ acting on piano gestures, and $FP_V$ acting on violin gestures). There is similarity between these gestures, because the articulations and the dynamics are the same. However, if we choose a {\em staccato} gesture for the violin, there is similarity only for the dynamic, but not for the articulation. The diagram involving non-similarity is not commutative. All transformations involved are homotopic, except the one required to pass from {\em legato} to {\em staccato} gestures. We can argue the Heuristic Conjecture \ref{conjecture_1}, connected with Definition \ref{theorem similarity}---its physical equivalent of Condition $2$ of the Definition.
\begin{conjecture}\label{conjecture_1}
Heuristic conjecture. Two gestures, based on the same skeleton, are {\em similar}\footnote{Two gestures are formally called {\em similar} if they satisfy the conditions of Definition \ref{theorem similarity}.} if and only if they can be connected via a transformation:
\begin{enumerate}
\item that homotopically transforms a gesture into the other,
\item and that also leads to similar changes in their respective acoustical spectra.
\end{enumerate}
Homotopy is a necessary, but not sufficient, condition to get similar gestures.
\end{conjecture}
Changes of loudness in orchestral playing lead also to changes in timbre, and so to changes in the spectra.
The tridimensional diagram of Figure \ref{3d_diagram} shows how we can compare gestures and their acoustical results inside the same diagram. In this way, the analysis of the spectrograms can be compared with the analysis of the gestural curves in space and time. Moreover, we may define a {\em degree of similarity} considering a parameter $\lambda$ that is equal to $1$ for two identical gestures, and decreases to zero for gestures more and more different.  

\begin{figure}\label{non_sim}
\centerline{\includegraphics[width=10cm]{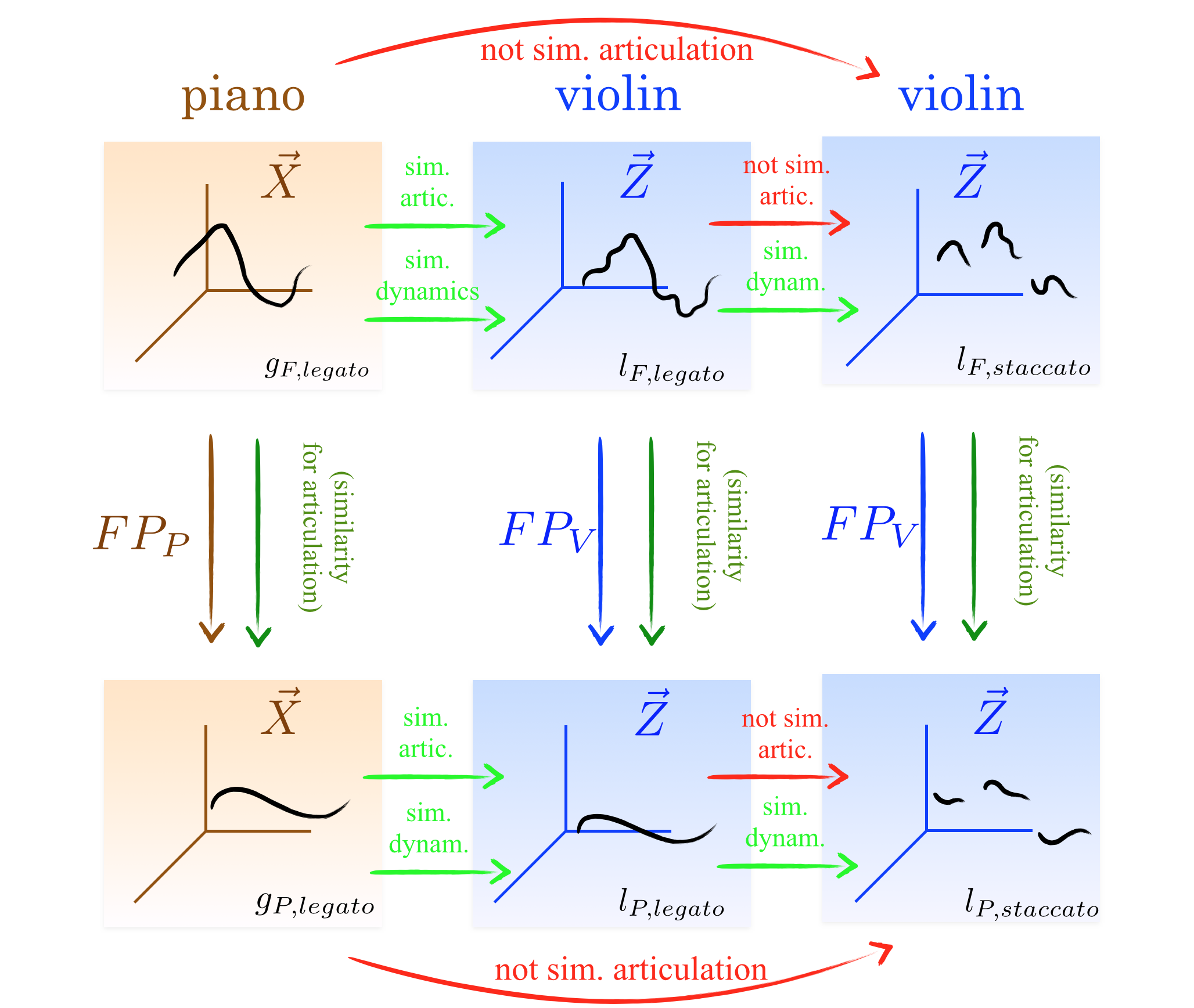}}
\caption{\label{non_sim} \footnotesize
Comparison between {\em forte} gestures on the piano and on the violin, {\em piano} gestures on both instruments, as well {\em legato} and {\em non legato}, as well transitions from forte to piano ({\em diminuendo}), and from {\em legato} to {\em staccato}.}
\end{figure}

\begin{figure}\label{3d_diagram}
\centerline{\includegraphics[width=16cm]{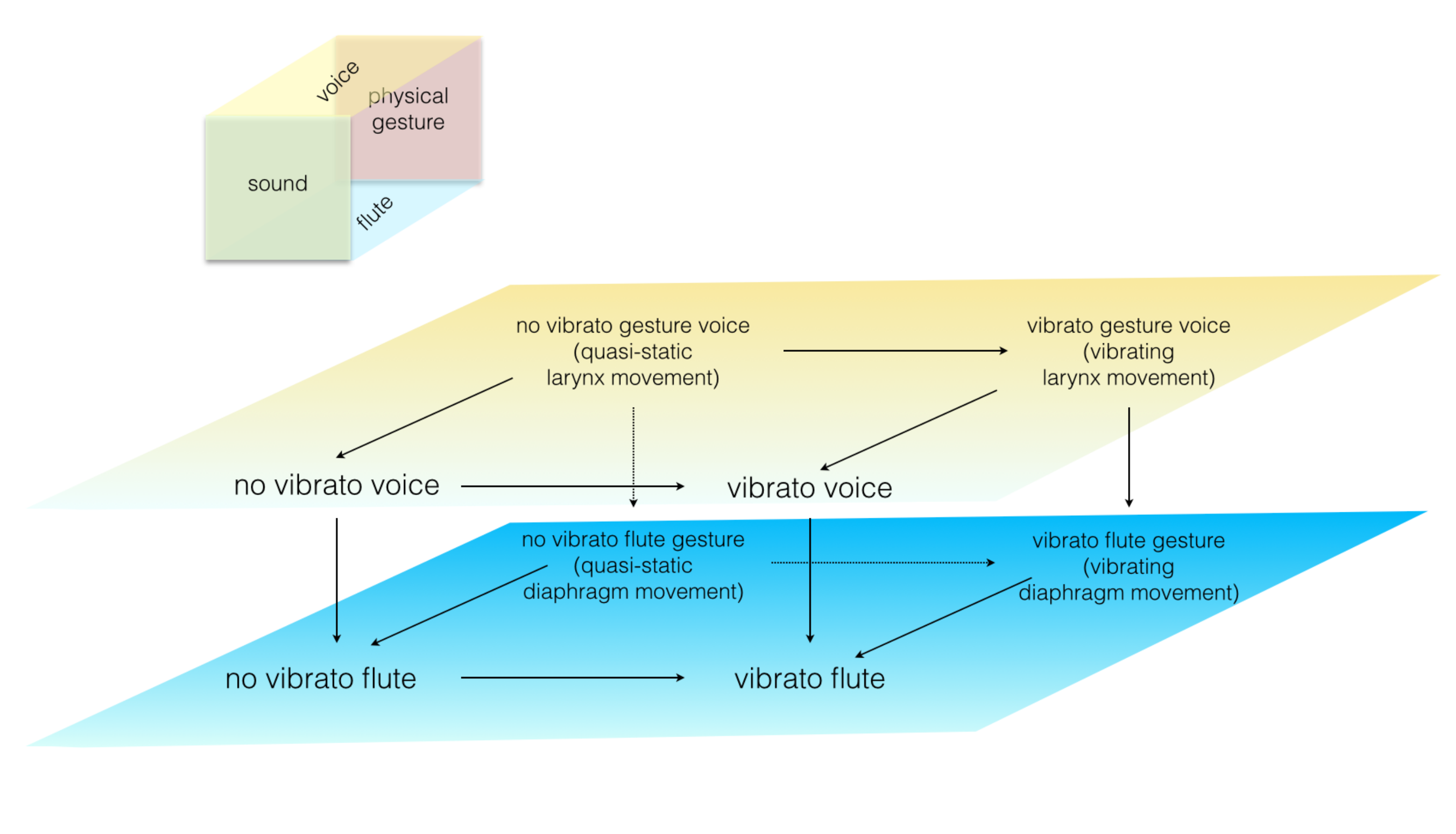}}
\caption{\label{3d_diagram} \footnotesize Tridimensional diagram that relates physical gestures with with their acoustical results \citep{phd_mannone}, in the case of the voice versus flute, and {\em vibrato} versus {\em non vibrato}. This kind of diagrams is not always commutative: it is not commutative, for example, when a sound  cannot be made by any human performer doing a musical gesture on an acoustic instrument.}
\end{figure}

\subsection{The category $\nabla$}\label{nabla}

We describe now the $\nabla$ category \citep{global}, used in Section \ref{gesti parametrizzati}. The notation $\nabla$ denotes here an internal category in $Top$, generalizing the unitary interval $I=[0,1]$ in the two-dimension plane. To define $\nabla$, let us think of the triangular region of space above the line $x=y$, consisting of points $(x,y)$ with $0\leq x<1$ and $0\leq y<1$.
The $\nabla$ category has:
\begin{enumerate}
\item as objects, the points belonging to the diagonal, i.e. their pairs of coordinates that verify $x=y$, and
\item as morphisms, the points inside the upper triangle, i.e. their pairs of coordinates verifying $x\leq y$, the pairs $(x,y)$ with $x\leq y$ of real numbers $x,y\in I=[0,1]$.
\end{enumerate}
The choice of the symbol $\nabla$ is due to such a triangular form.
We can use the symbol $\nabla$ to distinguish the set of points in the space $X$ and the set of arrows $X^{\nabla}$ having points in $X$, also indicated with $\vec X$.\footnote{Curves in $\vec X$ are mapped into points of $X$, for example with functions taking the tail and the head of an arrow. We can interpret these functions as directed graphs of the set of arrows and points.
In fact, we indicate as $\vec X$ a directed graph, identified by the curves $c:\nabla\rightarrow X$, whose the $head$ and $tail$ functions, projecting to $X$ the final and initial points of continuous curves in $\vec X$. This action of projecting a curve into points can be generalized not only for final and initial points, but also for every other point belonging to the curve.}


\subsection{Parametrized gestures}\label{gesti parametrizzati}

\begin{figure}
\centerline{\includegraphics[width=4.5cm]{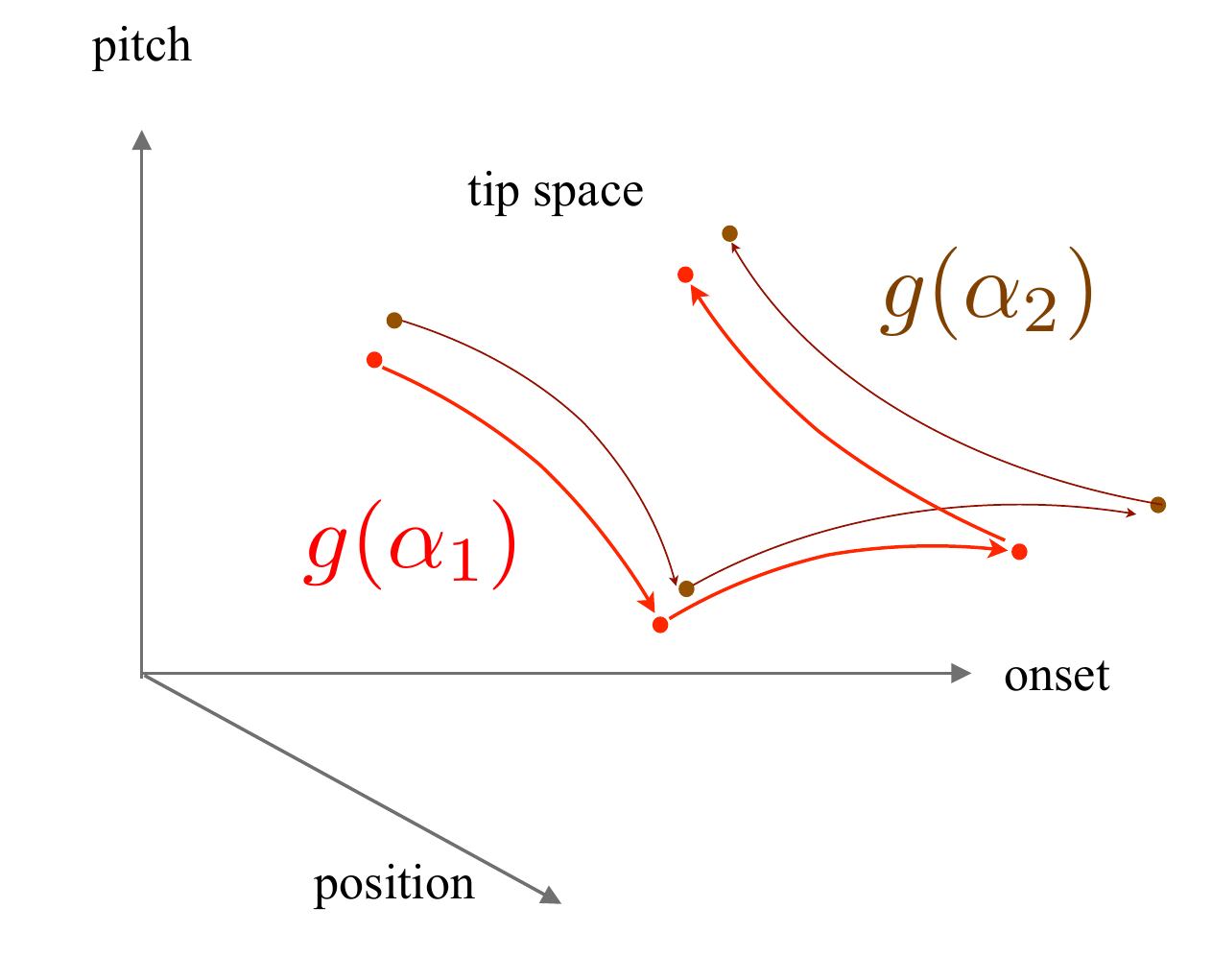}}
\caption{\label{parameter} \footnotesize Two slightly different parametrized gestures, $g(\alpha_1)$ and $g(\alpha_2)$, corresponding to a bit different choice of parameters $\alpha_1$ and $\alpha_2$. The {\em tip space} is the space of the positions of a finger's tip of the hand of the pianist.  A tip space indicated as position-onset-pitch has been used in \citep{mazzola_andreatta}.
}
\end{figure}

The concept of {\em parametrized gestures} is introduced
to formally take into account the influence of physical and physiological parameters on gestures.
\footnote{Parametrized gestures have been used to mathematically describe the mechanism of voice in singing \citep[Chapter 37]{ToM_III}. The dimension of the hand of the pianist, the anatomy of the vocal tract for the singer, and the values of diaphragm pressure and position of the larynx are examples of parameter choices.}
For each choice $\alpha$ of parameters in a parameter space $A$, we have a corresponding gesture $g(\alpha)$ in $\Delta\vec@ X$.
Each choice of parameters implies a different ``embodiment'' of the given skeleton $\Delta$ of Figure \ref{gesture}. Moreover, we require that the mapping such as $g:A\rightarrow\Delta\vec@X$ is continuous, and we define the set of these $A$-parametrized curves\footnote{We can define an {\em $A$-parametrized gesture} as a continuous function from the Cartesian product of the parameter space $\nabla$ with the category $A$,
with values in the topological space $X$, that is, $A\times\nabla\rightarrow X$. According to what is called the {\em currying theorem} in informatics, for a category $C$ we have $C(X\times Y,Z)\tilde\rightarrow C(X,Z^Y)$,
where $Z^Y$ are the curves from $Y$ to $Z$, and we can write that 
$q:A\times\nabla\rightarrow X\sim q:A\rightarrow\nabla@X$. As described in \citep{global}, an {\em $A$-addressed gesture with skeleton digraph $\Delta$ and body $X$} is a digraph morphism $g:\Delta \rightarrow A@\vec X$ into the spatial digraph of $A@\vec X$.} as $\Delta\vec@_AX$. Figure \ref{parameter} shows two slightly different gestures, corresponding to a slightly different choice of parameters $\alpha_1$ and $\alpha_2$.
More formally, let $\nabla$ be the topological category of curve parameters, and $A$ a topological category. For every choice of parameters, we can define a curve with values in the topological category $X$, given by
$q:A\rightarrow\nabla@X.$
This is equivalent to $q:A\times\nabla\rightarrow X$.
This means that, for $\alpha\in A$, we have a curve $C(\alpha):\nabla\rightarrow X$.
Let us consider the unit interval $I=[0,1]$ and a topological category $A$. 
Going along $I$ with a parameter $t$ ranging from $0$ to $1$, we can build a parametric curve having $s$ as parameter (for example: $s$ for the variable of $(\sin(s),\cos(s))$). In the rectangle $A\times I$, we have a set of values of $s$ ranging from $0$ to $1$ {\em for each choice} of parameter $\alpha\in A$. A straight line in $A\times I$ is {\em not yet} a skeleton: it can instead constitute one of its parts. These values can be mapped into a curve in the topological space $X$,
exactly as it happens in the definition of gesture, from a skeleton of a digraph to a curve (system of curves) in a topological space.

A {\em morphism} $f:g\rightarrow h$ of two {\em addressed gestures}  $g:\Gamma\rightarrow A@\vec X$ and
$h:\Delta\rightarrow B@\vec Y$, is a triple
$f = (t:\Gamma\rightarrow\Delta,\,a:B\rightarrow A,m:X\rightarrow Y)$, consisting of a digraph morphism $t$, an address change $a$, and a continuous functor of topological categories $m$, such that 
$\Gamma@_aY\circ\Gamma@_Am(g) = t@_BY(h)$.
Finally, we have that an $A$-addressed gesture with skeleton digraph $\Gamma$ and body $X$ is a digraph morphism $g:\Gamma\rightarrow A@\vec X$,
and we denote this set of gestures by $\Gamma@_A\vec X$.
We can finally write the relations of equations \ref{param1}, \ref{param2}, and \ref{param3} \citep{global}.
\begin{equation}\label{param1}
\Gamma@_aX:\Gamma@_A\vec X\rightarrow\Gamma@_B\vec X
\end{equation}
\begin{equation}\label{param2}
\Gamma@_Am:\Gamma@_A\vec X\rightarrow\Gamma@_A\vec Y \end{equation}
\begin{equation}\label{param3}
t@_AX:\Gamma@_A\vec X\rightarrow\Delta@_A\vec X.
\end{equation}

\subsubsection{Parametrized hypergestures}\label{hypergesture_par}

The concept of parametrization can also be applied to hypergestures. Let us indicate with $\Delta@X$ the set of all the curves from $\Delta$ to $X$. Let $g$ be a $B$-parametrized gesture,
going from the skeleton $\Delta$ to the system of curves in the space $X^B$. Briefly, $g\in\Delta\vec@_BX$, where we consider the new set having the structure of topological space.
In such a space, the gesture $g$ is a point. As known, the morphisms in this new space, connecting points (gestures), are hypergestures. If we have a skeleton $\Gamma$, mapped via the hypergesture $l$ into a system of $A$-parametrized curves connecting points $g$ and $h$ (both gestures $g,\,h\in\Delta\vec@_BX$), we say that $l$ is an $A$-parametrized hypergesture, and we can write $l\in\Gamma@_A(\Delta\vec@_BX)$.
If also this space has the structure of a topological space, we can write $\Gamma\vec@_A(\Delta\vec@_BX)$.



\subsection{Musical forces}\label{musical forces}

In the Introduction, we cited world-sheets. They can be framed in the topic of {\em musical forces} that shape their surface. World-sheets are obtained from the minimization of the action for a Lagrangian including complex time\footnote{In \citep{ToM_III} the {\em complex time} is introduced to extend the formalism of Physics to the mental reality of the score, where music flows in the imaginary time. The reality of the physical music performance flows in the real time \citep{ primas}.} \citep{mcm15, ToM_III}. The shape of the world-sheet is determined by an {\em (artistic) potential}, a force field.\footnote{The question of such an {\em artistic} potential is currently under research. It has an analogy at the level of notes with the performance operators \citep{gm:TOM, perfo}, both mathematically and conceptually.}
Such a surface creates a tight relation between symbolic and physical gesture, representing ``artistic'' forces that transform thought into artistic expression. The concept of force is a classic metaphor in musical creativity \citep{larson, gm:TOM, creativity, schoenberg, graben}. The concept of such a potential is used in numerical calculations, and used to mean the characteristic combination of parameters that distinguish a gesture from another \citep{global} but it has not been precisely defined yet.
The role of the force field may be related to the non-trivial conversion of the score's mental reality to the performance's physical reality. This process depends upon several elements such as technical difficulties, training, technical ability/skills of the performer, and expressivity of the piece.

\section*{Disclosure statement}
The author has no conflict of interest.

\section*{Acknowledgements}

I am grateful to Guerino Mazzola for the fruitful discussions about the relations between art and science in the framework of category theory.
I am also grateful to Peter beim Graben, for the reading of my manuscript, and for the insightful suggestions between physics, mathematics and semiotics. Thank you to Luca Nobile for reading and suggesting related studies in the field of crossmodal correspondences and neuroscience. I also would like to give a special thank to Giuseppe Metere for the mathematical insights. Finally, I am grateful to Marco Betta for musical advice, to Jason Yust for careful work as editor, and to Thomas Fiore for his careful corrections to the galley proofs. I also thank the referees.

\bibliographystyle{tMAM}
\bibliography{MySubmissionBibTexDatabase}

\begin{thebibliography}{57}
\newcommand{\enquote}[1]{``#1''}
\providecommand{\natexlab}[1]{#1}
\providecommand{\url}[1]{\normalfont{#1}}
\providecommand{\urlprefix}{ }
\expandafter\ifx\csname urlstyle\endcsname\relax
  \providecommand{\doi}[1]{doi:\discretionary{}{}{}#1}\else
  \providecommand{\doi}{doi:\discretionary{}{}{}\begingroup
  \urlstyle{rm}\Url}\fi

\bibitem[Alunni(2004)]{alunni}
Alunni, Charles. 2004.
\newblock \emph{Penser par le diagramme---de Gilles Deleuze {\`a} Gilles
  Ch{\^a}telet}, chap. Diagrammes et Cat{\'e}gories comme prol{\'e}gom{\`e}nes
  {\`a} la question : Qu'est-ce que s'orienter diagrammatiquement dans la
  pens{\'e}e ?
\newblock Presse Universitaire de Vincennes.

\bibitem[Antoniadis and Bevilacqua(2016)]{greek_2}
Antoniadis, Pavlos, and Fr{\'e}d{\'e}ric Bevilacqua. 2016.
\newblock ``Processing of symbolic music notation via multimodal performance
  data: Brian Ferneyhough's Lemma-Icon-Epigram for solo piano, phase 1.''
  \emph{Proceedings of TENOR 2016}
  \urlprefix\url{http://tenor2016.tenor-conference.org/papers/18_Antoniadis_tenor2016.pdf}.

\bibitem[Bachrat{\'a}(2011)]{bachrata}
Bachrat{\'a}, Petra. 2011.
\newblock ``Mod{\`e}les musicaux interactifs bas{\'e}s sur le geste pour
  l'analyse et la composition de musique mixte.'' \emph{Revue Francophone
  d'Informatique Musicale, Maison des Science de l'Homme, Paris Nord}  (1).
\newblock \urlprefix\url{http://revues.mshparisnord.org/rfim/index.php?id=123}.

\bibitem[Barbieri and Gentilucci(2009)]{gent1}
Barbieri, Buonocore Antonio Dalla Volta~Riccardo, Filippo, and Maurizio
  Gentilucci. 2009.
\newblock ``How symbolic gestures and words interact with each other.''
  \emph{Brain Lang} 110 (1).

\bibitem[beim Graben and Blutner(2016)]{graben}
beim Graben, Peter, and Reinhard Blutner. 2016.
\newblock ``Toward a Gauge Theory of Musical Forces.'' In \emph{Proceedings of
  the Quantum Interaction Conference 2016,} Lecture Notes in Computer Science.

\bibitem[beim Graben and Potthast(2009)]{graben2}
beim Graben, Peter, and Roland Potthast. 2009.
\newblock ``Inverse problems in dynamic cognitive modeling.'' \emph{Chaos} 19
  (1): 015103.

\bibitem[Berlin(2006)]{berlin}
Berlin, Brent. 2006.
\newblock ``The First Congress of Ethnozoological Nomenclature.'' \emph{Journal
  of the Royal Anthropological Institute} 12 (s1): S23--S44.

\bibitem[Bevilacqua et~al.(2011)Bevilacqua, Schell, Rasamimanana, Zamborlin,
  and Gu{\'e}dy]{bevilacqua}
Bevilacqua, Fr{\'e}d{\'e}ric, Norbert Schell, Nicolas Rasamimanana, Bruno
  Zamborlin, and Fabrice Gu{\'e}dy. 2011.
\newblock \emph{Musical Robots and Interactive Multimodal Systems}, chap.
  Online Gesture Analysis and Control of Audio Processing, 127--142.
\newblock Springer-Verlag, Basel.

\bibitem[Bevilacqua et~al.(2009)Bevilacqua, Zamborlin, Sypniewski, Schnell,
  Gu{\'e}dy, and Rasamimanana]{follower}
Bevilacqua, Fr{\'e}d{\'e}ric, Bruno Zamborlin, Anthony Sypniewski, Norbert
  Schnell, Fabrice Gu{\'e}dy, and Nicolas Rasamimanana. 2009.
\newblock ``{Continuous Realtime Gesture Following and Recognition}.'' In
  \emph{{Gesture in Embodied Communication and Human-Computer Interaction},}
  edited by Stefan Kopp and Ipke Wachsmuth, {Lecture Notes in Computer
  Science}. Springer.

\bibitem[Buteau and Mazzola(2000)]{mazzola_motif}
Buteau, Chantal, and Guerino Mazzola. 2000.
\newblock ``From Contour Similarity to Motivic Topology.'' \emph{Musicae
  Scientiae} VI (2).

\bibitem[Cont, Dubnov, and Assayag(2011)]{cont}
Cont, Arshia, Shlomo Dubnov, and G{\'e}rard Assayag. 2011.
\newblock ``On the information geometry of audio streams with applications to
  similarity computing.'' \emph{Audio, Speech, and Language Processing, IEEE
  Transactions on} 4: 837--846.

\bibitem[Fran{\c c}oise, Caramiaux, and Bevilacqua(2012)]{fran}
Fran{\c c}oise, Jules, Baptiste Caramiaux, and Fr{\'e}d{\'e}ric Bevilacqua.
  2012.
\newblock ``A Hierarchical Approach for the Design of Gesture to Sound
  Mappings.'' In \emph{Proceedings of the 9th Sound and Music Conference
  (SMC),} Copenhagen.

\bibitem[Gentilucci and Gianelli(2008)]{gent4}
Gentilucci, Dalla Volta~Riccardo, Maurizio, and Claudia Gianelli. 2008.
\newblock ``When the hand speak.'' \emph{The Journal of Physiology} 102:
  21--30.

\bibitem[Gentilucci and Bernardis(2006)]{gent3}
Gentilucci, Maurizio, and Paolo Bernardis. 2006.
\newblock ``Speech and gesture share the same communication system.''
  \emph{Neurospychologia} 44 (178-190).

\bibitem[Gentilucci and Bernardis(2007)]{gent5}
Gentilucci, Maurizio, and Paolo Bernardis. 2007.
\newblock ``Imitation during phoneme production.'' \emph{Neuropsychologia} 45
  (3): 608--615.

\bibitem[Gentilucci and Dalla~Volta(2008)]{gent2}
Gentilucci, Maurizio, and Riccardo Dalla~Volta. 2008.
\newblock ``Spoken language and arm gesture are controlled by the same motor
  control system.'' \emph{The Quarterly Journal of Experimental Psychology} 61
  (6): 944--957.

\bibitem[God{\o}y and Leman(2009)]{godoy_leman}
God{\o}y, Rolf, and Marc Leman, eds. 2009.
\newblock \emph{Musical Gestures: Sound, Movement and Meaning}.
\newblock Routledge, New York.

\bibitem[Henrotte(1992)]{semiotics}
Henrotte, Gayle. 1992.
\newblock ``Music and Gesture: A Semiotic Inquiry.'' \emph{American Journal of
  Semiotics} 9 (4): 103--113.

\bibitem[Iazzetta(2000)]{semiotics2}
Iazzetta, Fernando. 2000.
\newblock \emph{{Trends in Gestural Control of Music}}, chap. {Meaning in
  Musical Gesture}.
\newblock Ircam, Centre Pompidou.

\bibitem[Joyal(2008)]{Joyal}
Joyal, Andr{\'e}. 2008.
\newblock ``{The theory of quasi-categories and its applications}.''
  \emph{Centre de Recerca Matematica Barcelona} 45, II.
\newblock
  \urlprefix\url{http://mat.uab.cat/~kock/crm/hocat/advanced-course/Quadern45-2.pdf}.

\bibitem[Knees and Schedl(2016)]{new_book_sim}
Knees, Peter, and Markus Schedl. 2016.
\newblock \emph{Music Similarity and Retrieval}.
\newblock Springer, Heidelberg.

\bibitem[Kosko(1994)]{kosko}
Kosko, Bart. 1994.
\newblock \emph{Fuzzy Thinking: The New Science of Fuzzy Logic}.
\newblock Hyperion, Westport.

\bibitem[Kubovy and Schutz(2010)]{kubovy}
Kubovy, Michael, and Michael Schutz. 2010.
\newblock ``{Audio-visual objects}.'' \emph{Review of Philosophy and
  Psychology} 1 (1): 41--61.

\bibitem[Larson(2012)]{larson}
Larson, Steve. 2012.
\newblock \emph{Musical Forces}.
\newblock Musical Meaning and Interpretation. Indiana University Press,
  Bloomington.

\bibitem[Lawvere and Schanuel(2009)]{lawvere}
Lawvere, William, and Stephen Schanuel. 2009.
\newblock \emph{{Conceptual Mathematics. A first introduction to categories}}.
\newblock Cambridge University Press.

\bibitem[Leinster(1998)]{cat1}
Leinster, Tom. 1998.
\newblock ``Basic Bicategories.'' \emph{ArXiV}
  \urlprefix\url{https://arxiv.org/pdf/math/9810017v1.pdf}.

\bibitem[Lurie(2008)]{infinite}
Lurie, Jacob. 2008.
\newblock ``What is an $\infty$-Category?.'' \emph{Notices of the AMS} 55 (8):
  949--950.

\bibitem[Mac~Lane(1971)]{macLane}
Mac~Lane, Saunders. 1971.
\newblock \emph{Categories for the Working Mathematician}.
\newblock Springer, Heidelberg.

\bibitem[Mannone(2011)]{music_image_book}
Mannone, Maria. 2011.
\newblock \emph{Dalla Musica all'Immagine, dall'Immagine alla Musica. Relazioni
  matematiche fra composizione musicale e arte figurativa (From Music to
  Images, from Images to Music. Mathematical relations between musical
  composition and figurative art)}.
\newblock Palermo: Edizioni Compostampa.

\bibitem[Mannone(2017)]{phd_mannone}
Mannone, Maria. 2017.
\newblock ``{Musical Gestures between Scores and Acoustics: A Creative
  Application to Orchestra}.'' Ph.D. thesis, University of Minnesota.

\bibitem[Mannone and Compagno(2014)]{mannone_mus_nonMarkov}
Mannone, Maria, and Giuseppe Compagno. 2014.
\newblock ``Characterization of the degree of Musical non-Markovianity.''
  \emph{ArXiV} \urlprefix\url{http://arxiv.org/pdf/1306.0229v2.pdf}.

\bibitem[Mannone, Lo~Franco, and Compagno(2013)]{mannone_physics}
Mannone, Maria, Rosario Lo~Franco, and Giuseppe Compagno. 2013.
\newblock ``Comparison of non-Markovianity criteria in a qubit system under
  random external fields.'' \emph{Physica Scripta} T153: 014047.

\bibitem[Mannone and Mazzola(2015)]{mcm15}
Mannone, Maria, and Guerino Mazzola. 2015.
\newblock ``{Hypergestures in Complex Time: Creative Performance Between
  Symbolic and Physical Reality}.'' In \emph{Proceedings of the MCM 2015
  Conference,}   edited by T.~Collins et~al., 137--148. Springer, Heidelberg.

\bibitem[Mayerthaler(1987)]{mayer}
Mayerthaler, Willi. 1987.
\newblock ``System-independent morphological naturalness.'' In \emph{Leitmotifs
  in Natural Morphology,}  Vol.~10 of \emph{Studies in Language Companion
  Series}  edited by W.~U. Dressler, W.~Mayerthaler, O.~Panagl, and W.~U.
  Wurzel, 25--58. Benjamins.

\bibitem[Mazzola(2011)]{perfo}
Mazzola, Guerino. 2011.
\newblock \emph{Musical Performance}.
\newblock Computational Music Science. Springer, Heidelberg.

\bibitem[Mazzola and Andreatta(2007)]{mazzola_andreatta}
Mazzola, Guerino, and Moreno Andreatta. 2007.
\newblock ``{Diagrams, Gestures and Formulae in Music}.'' \emph{Journal of
  Mathematics and Music} 1 (1): 23--46.

\bibitem[Mazzola et~al.(To appear 2018)Mazzola, Guitart, Ho, Lubet, Mannone,
  Rahaim, and Thalmann]{ToM_III}
Mazzola, Guerino, Ren{\'e} Guitart, Jocelyn Ho, Alex Lubet, Maria Mannone, Matt
  Rahaim, and Florian Thalmann. To appear 2018.
\newblock \emph{The Topos of Music III: Gestures}.
\newblock Computational Music Science. Springer, Heidelberg.

\bibitem[Mazzola and Mannone(2016)]{global}
Mazzola, Guerino, and Maria Mannone. 2016.
\newblock ``Global Functorial Hypergestures over General Skeleta for Musical
  Performance.'' \emph{Journal of Mathematics and Music}
  \urlprefix\url{http://dx.doi.org/10.1080/17459737.2016.1195456}.

\bibitem[Mazzola, Park, and Thalmann(2011)]{creativity}
Mazzola, Guerino, Joomi Park, and Florian Thalmann. 2011.
\newblock \emph{Musical Creativity}.
\newblock Computational Music Science. Springer, Heidelberg.

\bibitem[Mazzola(2002)]{gm:TOM}
Mazzola, Guerino et~al. 2002.
\newblock \emph{The Topos of Music: Geometric Logic of Concepts, Theory, and
  Performance}.
\newblock Birkh\"auser Verlag, Basel.

\bibitem[Nobile(2013)]{nobile}
Nobile, Luca. 2013.
\newblock \emph{Iconicity in Language and Literature 10 : Semblance and
  Signification}, chap. Words in the mirror: Analysing the sensorimotor
  interface between phonetics and semantics in Italian, 101--131.
\newblock Amsterdam/Philadelphia: John Benjamins.

\bibitem[Nobile(2015)]{nobile2}
Nobile, Luca. 2015.
\newblock ``Phonemes as images: An experimental inquiry into shape-sound
  symbolism applied to the distinctive features of French.'' In
  \emph{Iconicity: East meets West,}   edited by Herlofsky William
  Shinohara~Kazuko Hiraga, Masako and Kimi Akita, 71--91.

\bibitem[Ohala(1982)]{ohala}
Ohala, John. 1982.
\newblock ``The voice of dominance.'' \emph{The Journal of the Acoustical
  Society of America} 72 (S66).

\bibitem[Primas(2007)]{primas}
Primas, Hans. 2007.
\newblock ``Non-Boolean descriptions for mind-matter problems.'' \emph{Mind and
  Matter} 5 (1): 7--44.

\bibitem[Rinman et~al.(2003)Rinman, Friberg, Kjellmo, Camurri, Cirotteau, Dahl,
  Mazzarino, Bendikesen, and McCarthy]{friberg}
Rinman, Marie-Louise, Anders Friberg, Ivar Kjellmo, Antonio Camurri, Damien
  Cirotteau, Sofia Dahl, Barbara Mazzarino, Bendik Bendikesen, and Hugh
  McCarthy. 2003.
\newblock ``EPS - An Interactive Collaborative Game Using Non-Verbal
  Communication.'' In \emph{Proceedings of the Stockholm Music Acoustics
  Conference,} Stockholm.

\bibitem[Roffler and Butler(1968)]{roffler}
Roffler, Suzanne, and Robert Butler. 1968.
\newblock ``Factors that influence the localization of sound in the vertical
  plane.'' \emph{The Journal of the Acoustical Society of America} 43:
  1255--1259.

\bibitem[Schoenberg(1983 edition)]{schoenberg}
Schoenberg, A. 1983 edition.
\newblock \emph{Theory of Harmony}.
\newblock University of California Press.

\bibitem[Schutz and Lipscomb(2007)]{lipscomb}
Schutz, Michael, and Scott Lipscomb. 2007.
\newblock ``Hearing gestures, seeing music: Vision influences perceived tone
  duration.'' \emph{Perception} 36: 888--897.

\bibitem[Sessions and Cone(1979)]{sessions2}
Sessions, Roger, and Edward~T. Cone. 1979.
\newblock \emph{Roger Sessions on Music, Collected Essays}.
\newblock Princeton Legacy Library.

\bibitem[Spence(2011)]{spence}
Spence, Charles. 2011.
\newblock ``Crossmodal correspondences: A tutorial review.'' \emph{Attention,
  Perception and Psychophysics} 73: 971--995.

\bibitem[Spivak(2013)]{spivak}
Spivak, David. 2013.
\newblock \emph{{Category Theory for Scientists}}.
\newblock MIT.

\bibitem[Uznadze(1924)]{uznadze}
Uznadze, Dimitri. 1924.
\newblock ``Ein experimenteller Beitrag zum Problem der psychologischen
  Grundlagen der Namengebung (An experimental contribution to the problem of
  the psychological foundations of naming).'' \emph{Psychologische Forschung} 5
  (1-2): 24--43.

\bibitem[Visi, Schramm, and Miranda(2014)]{visi_miranda_2}
Visi, Federico, Rodrigo Schramm, and Eduardo Miranda. 2014.
\newblock ``Gesture in performance with traditional musical instruments and
  electronics: Use of embodied music cognition and multimodal motion capture to
  design gestural mapping strategies.'' In \emph{Moco '14: proceedings of the
  2014 international workshop on movement and computing,} Paris.

\bibitem[Wieland and Uhde(2002)]{wieland}
Wieland, Renate, and J{\"u}rgen Uhde. 2002.
\newblock \emph{{Forschendes {\"U}ben. Wege instrumentalen Lernens. {\"U}ber
  den Interpreten und den K{\"o}rper als Instrument der Musik}}.
\newblock B{\"a}renreiter, Neum{\"u}nster.

\bibitem[Winter(2004)]{fuzzy1}
Winter, Michael. 2004.
\newblock ``Goguen Categories.'' \emph{Journal on Relational Methods in
  Computer Science} 1: 339--357.

\bibitem[Zbikowski(2002)]{zbikowski}
Zbikowski, Lawrence. 2002.
\newblock \emph{Conceptualizing Music}.
\newblock AMS Studies in Music. Oxford University Press.

\bibitem[Zwiebach(2004)]{zwiebach}
Zwiebach, Barton. 2004.
\newblock \emph{{A First Course in String Theory}}.
\newblock Cambridge University Press, Cambridge.

\end{thebibliography}

\addcontentsline{toc}{section}{References}

\end{document}